\documentclass[11pt]{article}
\input{amssym.def}
\input{amssym}
\usepackage[active]{srcltx}
\usepackage{color}

\usepackage{hyperref}
\hypersetup{
    colorlinks=true, 
    linktoc=all,     
    linkcolor=blue,  
}

\newenvironment{annotacia}{\centerline{\sc Abstract}\vspace{2mm}\narrower\narrower\sf}

\def\theequation{\thesection.\arabic{equation}}
\makeatletter\@addtoreset{equation}{section}\makeatother

\def\nn{\nonumber}\def\lb{\label}\def\nin{\noindent}
\def\be{\begin{equation}}\def\ee{\end{equation}}\def\ba{\begin{eqnarray}}\def\ea{\end{eqnarray}}
\def\tr{{\rm Tr}\,}
\def\Tr#1{{\rm Tr}_{\! R^{\mbox{\scriptsize$(#1)$}}}}
\def\TR#1#2{{\rm Tr}_{\! #2^{\mbox{\,\scriptsize$(#1)$}}}}
\def\str#1{\rule[#1mm]{0pt}{1mm}}

\def\rdet{{\rm{det}}}

\def\dM{\dot{M}}
\def\ddM{\ddot{M}}

\newcounter{theorem}\makeatletter
\@addtoreset{theorem}{section}\makeatother

\newtheorem{prop}[theorem]{Proposition}\newtheorem{rem}[theorem]{Remark}
\newtheorem{lem}[theorem]{Lemma}
\newtheorem{def-lem}[theorem]{Definition-Lemma}
\newtheorem{def-prop}[theorem]{Definition-Proposition}
\newtheorem{defin}[theorem]{Definition}\newtheorem{theor}[theorem]{Theorem}
\newtheorem{cor}[theorem]{Corollary}

\begin{document}
\title{ }
\begin{center}
{\Large \textbf{Reciprocal Relations for Orthogonal\\[3pt]
Quantum Matrices}}
	
\vspace{1cm} {\large \textbf{Oleg Ogievetsky$^{\diamond\,\dag}$ 
and Pavel Pyatov$^{\ast\,\star
}$}}
	
\vskip .8cm $^{\diamond}$Aix Marseille Universit\'{e}, Universit\'{e} de
Toulon, CNRS, \\ CPT UMR 7332, 13288, Marseille, France
	
\vskip .3cm $^{\dag}${I.E.Tamm Department of Theoretical Physics, 
P.N. Lebedev Physical Institute, Leninsky prospekt 53, 119991 Moscow, Russia}
		
\vskip .3cm $^{\ast}${National Research University "Higher School 
of Economics", 20 Myasnitskaya street, Moscow
101000, Russia}
	
\vskip .3cm $^{\star}${Bogoliubov Laboratory of Theoretical Physics, 
JINR, 141980 Dubna, Moscow region, Russia}
\end{center}

\begin{annotacia}\noindent 
For the family of the orthogonal quantum matrix algebras we investigate the structure of their  
characteristic subalgebras --- special commutative subalgebras, which for the 
subfamily of the reflection equation algebras appear to be central.
In \cite{OP-BMW} we described three generating sets of the characteristic subalgebras 
of the symplectic and orthogonal quantum matrix algebras. 
One of these --- the set of the 
elementary sums --- is finite. 
In the symplectic case the elementary sums are in general 
algebraically independent. On the contrary, 
in the orthogonal case the elementary sums turn out to be 
dependent. We obtain a set of quadratic reciprocal relations for these generators. 
Next, we resolve the reciprocal relations for the quantum orthogonal matrix algebra extended
by the inverse of the quantum matrix.
As an auxiliary result, we derive the commutation  
relations between the $q$-determinant of the
quantum orthogonal matrix and the generators of the 
quantum matrix algebra, that is, the components of the quantum matrix.
\end{annotacia}

\newpage
\tableofcontents
\bigskip\bigskip\bigskip
\section{Introduction}\lb{sec1}
A concept of Quantum Matrix (shortly, QM-)algebra unifies a  variety of 
algebras arising in the theory of quantum groups, among them most famous are 
families of the FRT-algebras \cite{FRT} and the Reflection Equation (RE-)algebras 
\cite{Ch,KS}. Roughly speaking QM-algebra is a unital associative algebra generated by  
the matrix components with the defining quadratic relations. These relations are expressed 
with the help of a Yang--Baxter representation 
of the braid group (for exact definitions see sections \ref{subsec2.2} and 
\ref{subsec4.1}).
Fruitfullness of the concept became appparent when it was realized that for these algebras 
one can universally define generalizations of the matrix power, the matrix trace and  
determinant, and, finally, derive a generalization of the fundamental theorem of the 
matrix algebra --- the Cayley-Hamilton theorem \cite{IOP1,IOP2}. In so doing, all the 
notions and results rely exclusively on the properties of the Yang--Baxter matrices appearing in 
the definiton of the QM-algebra.

So far, one usually considers the QM-algebras based on the 
representations of either  Iwahori--Hecke \cite{K,Iw}, or 
Birman--Murakami--Wenzl algebras \cite{BW,Mr}. In the first case we mean 
the so called $GL(m|n)$ type 
Yang--Baxter matrices (including the subseries of the type $GL(n)$); the second case decomposes
into the subtypes of symplectic $Sp(2n)$ Yang--Baxter matrices, or orthogonal $O(n)$  
Yang--Baxter matrices.\footnote{For examples of these classical series of the $R$-matrices the 
reader is referred  to \cite{FRT,I}} The more general case of orthosymplectic 
$OSp(m\vert 2n)$ Yang--Baxter matrices is less investigated.
 
We call the QM-algebras by the types of 
the $R$-matrices used in their construction.
The family of $GL(n)$ type QM-algebras is best studied. The 
structure results, including the analogue of the Cayley--Hamilton theorem, its generalization, 
Cayley--Hamilton--Newton theorem \cite{IOP}, and the thorough  
investigation of the special commutative subalgebra of the QM-algebra, 
called the {\em characteristic subalgebra}, were carried out in
\cite{IOP1,IOP2}. Many particular examples of the $GL(n)$ type QM-algebras 
were considered in details (see, e.g., \cite{EOW,NT,PS,GPS1,Zh,JW,OV}). The 
family of $GL(m|n)$ type QM-algebra was investigated in the general setting in 
\cite{GPS2,GPS3}. 

The program of investigation of the general orthogonal and symplectic types 
QM-algebras was initiated in the unpublished work 
\cite{OP}; our main objective was the generalization of the Cayley--Hamilton theorem for the quantum symplectic and orthogonal groups. Later, in \cite{OP-BMW} the characteristic 
subalgebra of these QM-algebras was described in terms of generators and relations. 
Various sets of generators were considered and relations among them were derived. 
In the paper \cite{OP-SpCH} on the base of obtained results we proved an analogue of the 
Cayley-Hamilton theorem for the symplectic QM-algebras.
Our next aim is a generalization of the Cayley-Hamilton theorem for the
orthogonal QM-algebras \cite{OP3}. In the orthogonal case the 
generators of the characteristic subalgebra are not algebraically independent and we 
need to find additional relations among them. These relations, called {\it reciprocal} 
relations will be necessary for the proper understanding of the fundamental notion of the 
eigenvalue of the
orthogonal quantum matrix. These reciprocal relations are derived in the present work.


\vskip .2cm
Before proceeding with considerations let us make a remark about the subfamily of 
RE-algebras. They are distinguished among the QM-algebras because
their characteristic subalgebras are central (for the proof see \cite{I}, proposition 5). 
Besides, they are closely related to the quantized universal enveloping algebras (see 
\cite{Sem} and, e.g., \cite{JW}, section 1.4). Accordingly, earlier results on characteristic 
identities for the 
quantized universal enveloping algebras \cite{GZB,MRS} and for their classical 
analogs (for a review see \cite{IWG}) agree with the results obtained for the 
RE-algebras. Structure results on the RE-algebras have found numerous applications. 
These results include but are not limited to, 
the quantization of the conjugacy classes 
\cite{DM,Md,AM1,AM2}, investigation of differential geometric constructions over 
quantum groups \cite{IP,P}, description of the diagonal reduction algebras \cite{KhO}, 
and the construction of the generalized Capelli identities \cite{GPeS,JLM,Z}.
\medskip

A word about notation. The main objects in this paper are defined with the help of a 
matrix solution of the Yang--Baxter equation which is an operator in the tensor square
$V\otimes V$ where $V$ is a vector space. Such solution is often called R-matrix;
however we shall call it a {\it Yang--Baxter} matrix to avoid word combinations 
like ``R-matrix R". Also, we shall call a braid group representation, naturally built on a 
Yang--Baxter
matrix, a {\it Yang--Baxter} representation.  \medskip

The main text of the paper is organized as follows.
\medskip

Section \ref{sec2} contains necessary information about Birman--Murakami--Wenzl 
(BMW, for short) algebras and their Yang--Baxter representations.

The BMW algebra ${\cal W}_{n}(q,\mu)$ is defined in subsection \ref{BMW-1}. Three 
sets of mutually orthogonal idempotents in ${\cal W}_{n}(q,\mu)$ --- the $q$-antisymmetrizers 
$a^{(n)}$, the $q$-symmetrizers $s^{(n)}$ and the contractors $c^{(n)}$ --- are 
introduced here. Yang--Baxter representations of these elements are used in section 
\ref{sec4} in the construction of  commutative subalgebras in the QM-algebras.

Yang--Baxter representations of the BMW algebras themselves, and the Yang--Baxter 
matrix  techniques including the notions of the R-trace and of the compatible Yang--Baxter matrix  
pair are introduced in subsections \ref{subsec2.2} and \ref{RF-pairs}. Formulas 
specific for the BMW type
Yang--Baxter matrices are collected here. 

\medskip
Orthogonal ($O(k)$ type) Yang--Baxter matrices are introduced in section \ref{sec3}. 
This particular family of the BMW type Yang--Baxter matrices
is used in the definition of the orthogonal  QM-algebras. In subsection \ref{subsec3.5} a 
series of the rank=1 projectors related to the orthogonal Yang--Baxter matrices are 
investigated. These projectors serve as a main technical tool in considerations in the 
next section.
\medskip

Section \ref{sec4} begins with a general definition the QM-algebra ${\cal M}(R,F)$ 
associated to a compatible pair $\{R,F\}$ of the strict skew invertible Yang--Baxter 
matrices. Imposing additional conditions on the Yang--Baxter matrix  $R$ one  
specializes to the QM-algebras of the BMW  and of the orthogonal types.  Note that in 
our definition the BMW type QM-algebras (and also their orthogonal and symplectic 
subfamilies) contain a nontrivial quadratic element --- the {\em contraction} $g$. 
Traditionally, in  consideration of examples of such algebras this element is set 
proportional to $1$. Explicitly, $g=\mu^2 1$ in the notations of (\ref{tau2})
(see, e.g., \cite{FRT}, definition 11, or \cite{KSch}, section 9.3.1). In our studies we 
keep the contraction free, not only for sake of generality, but mainly 
because fixing the value of $g$ makes our derivations less clear.

Next, in subsection \ref{subsec4.1b} we define special commutative subalgebra 
${\cal C}(R,F)\subset {\cal M}(R,F)$ --- the {\em characteristic subalgebra} --- and 
introduce three sets of its generating elements: the so-called power,  
elementary and complete sums. In the end of the subsection  we 
define a natural extension of the BMW type QM-algebra containing the inverse of the 
quantum matrix $M$. Material of this subsection is mainly based on \cite{OP-BMW}. 
\smallskip

Our main result is contained in subsection \ref{subsec4.11}. Here in theorem 
\ref{theorem4.9} we derive algebraic relations specific for the elementary sums
 in the orthogonal QM-algebras. These relations are called {\em 
reciprocal}. 

\smallskip
The subsection \ref{qdetorth} is devoted to the investigation of a distinguished element of 
the orthogonal type QM-algebra --- the {\em quantum determinant}. We obtain 
the commutation relations of the
quantum determinant  with the QM-algebra generators. These relations involve the
operator $O\in {\rm Aut}(V)$ associated with a compatible pair $\{R,F\}$, where $R$ is 
of the orthogonal type.

\smallskip
Two representations, $\alpha_+$ and $\alpha_-$, of a general QM-algebra are 
considered in the subsection \ref{repqma}. We calculate images of the contraction 
and of the quantum determinant in these representations. As discussed in subsection 
\ref{components} these data are useful for a resolution of the reciprocal relations.
In the end of the subsection \ref{repqma}, using the representation $\alpha_+$, we 
exhibit some properties of the operator $O$.

\smallskip
Our second main result is presented in the final subsection \ref{components}.  In theorem 
\ref{theorem4.23} we resolve the reciprocal relations for the extended orthogonal 
QM-algebras (that is, for the algebras extended by the inverse $g^{-1}$)
to the contraction $g$. 
The solution is sensitive to the parity of $k$ of the orthogonal $O(k)$ type QM-algebra.
For even $k$, $k=2\ell$, one has two different solutions associated with the two 
`components' of the type $O(2\ell)$ QM-algebra. These components are nothing else but the 
quantum analogues of the two connected components of the orthogonal matrix 
group. In the odd case, $k=2\ell -1$, the solution is unique. 

\vskip .2cm
The symbol $\blacktriangleright$ denotes ends of remarks as well as ends of proofs of auxiliary lemmas placed inside the proofs of theorems (whose ends are 
denoted by $\blacksquare$).

\section{Preliminaries}\lb{sec2}

In this section we recall necessary definitions and formulas concerning the 
Birman--Murakami--Wenzl algebras and their Yang--Baxter representations.
This material was discussed in detail in our preceding works 
\cite{OP-BMW,OP-SpCH}.
\smallskip
\subsection{BMW algebras}
\lb{BMW-1}
The {\em Birman--Murakami--Wenzl (BMW) algebra} ${\cal W}_{n}(q,\mu)$ 
\cite{BW,Mr} is a finite
dimensional quotient algebra of the group algebra ${\Bbb C}[{\cal B}_{n}]$ of the braid 
group ${\cal B}_{n}$. It can be defined in terms of generators
 $\{\sigma_i, \kappa_i\}_{i=1}^{n-1}$ and  relations depending on
two complex parameters $q$ and $\mu$:
\ba
\lb{braid group}
\sigma_i \sigma_{i+1} \sigma_i \, =\,  \sigma_{i+1} \sigma_i \sigma_{i+1},&&
\sigma_i \sigma_j \, =\,  \sigma_j \sigma_i\qquad   \forall\; i,j:\; |i-j|>1 ,
\\[2pt]
\lb{kappa}
\kappa_{i} = { \textstyle {(q1-\sigma_i)(q^{-1} 1 +\sigma_i)\over \mu (q-q^{-1})}} , 
&&
\sigma_i \kappa_i \, =\,  \kappa_i \sigma_i \,=\, \mu \kappa_i ,
\\[2pt]
\lb{bmw3} 
\kappa_i \sigma_{i \pm 1}^{\epsilon}  =\kappa_i \kappa_{i \pm 1}\sigma_i^{-\epsilon}\,
,&& \sigma_{i\pm 1}^{\epsilon} \kappa_{i}  =
\sigma_{i}^{-\epsilon} \kappa_{i\pm 1} \kappa_{i}\, ,\quad  \forall\,\epsilon \in\{1,-1\},
\\[2pt]
\lb{bmw5a} 
\kappa_i \kappa_{i\pm 1} \kappa_i = \kappa_i\, , &&
\kappa_i^2 = \eta\, \kappa_i\, ,
\quad
\mbox{where}\;\;
\eta:= {\textstyle {(q-\mu)(q^{-1}+\mu)\over\mu (q-q^{-1})}} .
\ea
Here $\epsilon$ and $\pi$ are signs:
$\epsilon =\pm 1$ and $\pi =\pm 1$.
The first line is the Artin's presentation of the braid group ${\cal B}_n$; the rest of
relations define the quotient algebra ${\cal W}_{n}(q,\mu)$ of the group algebra 
${\Bbb C}[{\cal B}_{n}]$.
The parameters $q$ and $\mu$  are taken in domains\footnote{For particular 
values  $\mu=\pm q^{i}$, $i\in {\Bbb Z}$, the limiting cases $q\rightarrow \pm 1$ to the 
Brauer algebra \cite{Br} can be consistently defined.}
\be
\lb{init-cond}
q\in{\Bbb C}\backslash \{0,\pm 1\},\quad \mu\in{\Bbb C}\backslash\{0,q,-q^{-1}\}, 
\ee
so that the elements $\kappa_i$ are well defined and
non-nilpotent. 

Inroduce two
sets of the so-called {\em baxterized elements} $\sigma^{\pm}_i$ \cite{Jon}:
\be
\lb{baxterization}
\sigma_i^{\pm}(x)\, :=\,  1\, +\, {x-1\over q-q^{-1}}\, \sigma_i\, +\,
{\mu(x-1)\over \mu\mp q^{\mp 1}  x}\, \kappa_i\, ,\quad 
x\in{\Bbb C}\setminus\{\pm q^{\pm 1} \mu\}.
\ee

Using the baxterized elements one defines two series of the idempotents in 
${\cal W}_{n}$: $a^{(i)}$ and $s^{(i)}$, $i=1,\dots, n$, called, respectively,
the {\em $q$-antisymmetrizers} and the {\em 
$q$-symmetrizers}\,\footnote{Coincidence of the two exressions 
for $q$-(anti)symmetrizers in (\ref{a^k}),(\ref{s^k}) and their further properties 
are explained in \cite{OP-BMW}, proposition 2.2.}
\ba
\lb{ind1}
a^{(1)} := 1\, , && s^{(1)}\, :=\, 1\, ,
\\[2pt]
\lb{a^k}
a^{(i+1)} \,\,:=\,\,{q^i\over (i+1)_q} a^{(i)}\, \sigma^{-}_i(q^{-2i})\,
a^{(i)} &
=&{q^i\over (i+1)_q} a^{(i)\uparrow 1}\, \sigma^{-}_1(q^{-2i})\,
a^{(i)\uparrow 1}\, ,
\\[2pt]
\lb{s^k}
s^{(i+1)}\,\, :=\,\,{q^{-i}\over (i+1)_q}\, s^{(i)}\, \sigma^{+}_i(q^{2i})\,
s^{(i)}  &
=&
{q^{-i}\over (i+1)_q}\, s^{(i)\uparrow 1}\,
\sigma^{+}_1(q^{2i})\, s^{(i)\uparrow 1}\, .
\ea
Here $i_q$ are the $q$-numbers, $i_q\, :=\, (q^i - q^{-i})/(q-q^{-1})$.
Symbol $x^{\uparrow 1}$ in these formulas stands for the image of the element 
$x\in {\cal W}_{n}$ (or, more generally, $x\in {\cal B}_n$) under the 
monomorphism ${\cal W}_{n}\hookrightarrow{\cal W}_{n+1}$ (respectively, 
${\cal B}_{n}\hookrightarrow{\cal B}_{n+1}$) defined on the generators as 
$(\sigma_i)^{\uparrow 1}= \sigma_{i+1}$~  $\forall i=1,\dots,k-1$.
Further on symbol $x^{\uparrow p}$ will denote $p$-times iteration of this 
map: $x^{\uparrow p}:= (x^{\uparrow (p-1)})^{\uparrow 1}$.

To avoid singularities in the definition of $a^{(i)}$ (respectively, 
$s^{(i)}$), $i=2,\dots ,n$,
we impose additional restrictions on the parameters of ${\cal W}_n(q,\mu)$:
\be
\lb{q-mu-restrict}
q^{2j}\neq\, 1, \quad \mu \neq -q^{-2j+3} \; (\mbox{resp.,}\
\mu \neq q^{2j-3})\, \; \forall\;j = 2,\dots , n\ .
\ee

In \cite{OP-BMW} yet another series of the idempotents was introduced:

\be
\lb{kappa-i}
 c^{(0)}\, :=\, 1, \qquad c^{(2)}\, :=\, \eta^{-1}\kappa_1\, , \quad\;\; c^{(2i+2)}\, :=\,
c^{(2i)\uparrow 1}\, \kappa_1 \kappa_{2i+1}\, c^{(2i)\uparrow 1}\;\; \forall\, i\geq 1 .
\ee
The element $c^{(2i)}$ is called the {\em $2i$-th order contractor}.
Alternatively, the element $c^{(2i)}$ can be defined by
\be
c^{(2i)} =
\eta^{-1} c^{(2i-2)\uparrow 1}\left(\kappa_{2i-1}\kappa_{2i-2}\dots\kappa_{i+1}\right)
\left(\kappa_1\kappa_2\dots\kappa_i\right)\;\;\forall\, i\geq 1 .
\lb{kappa-i2}
\ee
Equivalence of the two definitions follows from the relations (see \cite{OP-BMW})
\be
\lb{idemp-c1}
c^{(2n)} c^{(2i)\uparrow n-i} =c^{(2i)\uparrow n-i} c^{(2n)} = c^{(2n)}\quad
\forall\; i=1,2,\dots , n\, .
\ee
One also has
\be
\lb{idemp-c2}
c^{(2n)} \sigma_i = c^{(2n)} \sigma_{2n-i}\, , \qquad
\sigma_i c^{(2n)} = \sigma_{2n-i}\, c^{(2n)} \, \quad \forall\; i=1,2,\dots ,n-1\, .
\ee

\subsection{Yang--Baxter representations and R-trace}
\lb{subsec2.2}
We remind basic notations and definitions on the Yang--Baxter representations 
of the BMW algebras. A detailed exposition can be found in \cite{OP-BMW}, section 3. 

Let $V$ be a finite dimensional ${\Bbb C}$-linear space, 
$\dim V = \mbox{\sc n}$; ~$I\, :={\rm Id}_{_V}$.
An operator $R\in {\rm Aut}(V^{\otimes 2})$ is called 
an {\em Yang--Baxter matrix } (respectively, an {\em Yang--Baxter matrix  of a 
BMW type}) if the map defined on the set of the braid group (respectively, BMW 
algebra) generators 
\be 
\sigma_i \mapsto   R_i:= I^{\otimes(i-1)}\otimes R\otimes  I^{\otimes(n-i-1)
}\in {\rm Aut}(V^{\otimes n}), \quad 1\leq i\leq n-1,
\lb{localrep}
\ee
extends to a homomorphism
\be 
\rho_R:\, {\cal B}_n\; (\mbox{respectively,}\; {\cal W}_n)\rightarrow 
{\rm End}(V^{\otimes n})\ .
\lb{nequals3}
\ee

The homomorphism $\rho_R$ is called the {\em Yang--Baxter representation} of the 
group ${\cal B}_n$ (resp., the algebra ${\cal W}_n$).
Image of the braid relation (\ref{braid group}) in the representation $\rho_R$
\be 
R_1\, R_2\, R_1\, =\, R_2\, R_1\, R_2\, ,
\lb{YB}
\ee
is traditionally called the {\em braid relation}. This matrix equality is sufficient for the 
definition of the Yang--Baxter matrix  $R$ (that is, it is enough to take $n=3$ in 
(\ref{nequals3})).  If $R$ is the Yang--Baxter matrix, then so is $R^{-1}$.
\smallskip

An operator $R\in {\rm End}(V^{\otimes 2})$ is called {\it skew invertible} if there exists 
an operator
${\Psi_{\!R}}\in {\rm End}(V^{\otimes 2})$ such that
\be\tr_{(2)} R_{1} (\Psi_{\!R})_{2} =\tr_{(2)} (\Psi_{\!R})_{1} R_{2} = P .
\lb{s-inv}
\ee
Here $(\Psi_{\!R})_i$, $i=1,2$, are treated as elements of ${\rm End}(V^{\otimes 3})$, 
like in (\ref{localrep}), $(\Psi_{\!R})_1:=\Psi_{\!R}\otimes I$,
 $(\Psi_{\!R})_2:= I\otimes \Psi_{\!R}$.
The symbol ${\rm Tr}_{(i)}$ stands for the trace over an $i$-th component in the tensor 
product
of the spaces $V$.
The symbol $P\in {\rm Aut}(V^{\otimes 2})$ is reserved for the permutation operator: 
$P(u\otimes v)= v\otimes u \;\; \forall\, u, v\in V$. In the following, for an operator $Z$ in 
the tensor power of the space $V$, ${\rm Tr}_{(i_1,i_2,\dots ,i_{L})}$ stands for the 
trace over the components $i_1,i_2,\dots ,i_L$ in the tensor product.
\smallskip

With the skew invertible  Yang--Baxter matrix  $R$ one associates a pair of 
matrices
\be
\lb{C-D}
{(C_R)}_2:= \tr_{(1)} (\Psi_{\!R})_{1}\in {\rm End}(V) \,  \quad 
\mbox{and}\quad {(D_R)}_1:= \tr_{(2)} (\Psi_{\!R})_1\in {\rm End}(V)\, .
\ee
The meaning of the indices 1 and 2 in the symbols ${(D_R)}_1$ and ${(C_R)}_2$ is as 
follows. 
For an operator $Z\in{\rm End}(V)$ the symbol $Z_i$ stands for the operator  which 
acts as $Z$ in the $i$-th
component in the tensor power of the space $V$ and as the identity in all other 
components. In the sequel we shall often omit the brackets in the 
symbols ${(D_R)}_1$ and ${(C_R)}_2$ and write simply ${D_R}_1$ and ${C_R}_2$.

\medskip
The characteristic properties of operators $C$ and $D$ are, respectively,
\be
\lb{trC-D}
\tr_{(1)} {C_R}_1 R_{1} = {I}_2 \quad\mbox{and}\quad\tr_{(2)} {D_R}_2 R_{1} = {I}_1.
\ee
As shown in, e.g., \cite{O}, invertibility of (any of) the matrices $C_R$, $D_R$ is 
equivalent to the skew invertibility of the Yang--Baxter matrix  $R^{-1}$  and moreover
\be
\lb{ss-inv}
C_R\, D_{R^{-1}}\, =\, D_R\, C_{R^{-1}}\,  =\, I.
\ee
Skew invertible Yang--Baxter matrix  whose inverse matrix is also skew invertible is 
called {\em strict skew invertible}. 
For a Yang--Baxter matrix  $R$ of the BMW type, one has \cite{IOP3}
\be
\lb{Cinv-D}
C_R\, D_R\, =\, \mu^2\, I
\ee
and thus,  for the BMW type Yang--Baxter matrices skew invertibility implies strict skew 
invertibility.
\smallskip

With the matrix $D_R$ we associate a linear map on the space of 
$\mbox{\sc n}\!\times\! \mbox{\sc n}$ matrices whose entries belong to a 
$\Bbb{C}$-linear space $W$
$$
{\rm Tr\str{-1.3}}_R:\; {\rm End}(V)\otimes W\,\rightarrow \, W ,\;\; 
{\rm Tr\str{-1.3}}_R(X) ={\textstyle \sum_{i,j=
1}^{\mbox{\footnotesize\sc n}}}{(D_R)}_i^jX_j^i\, , \;\; 
X\in{\rm End}(V)\otimes W\,,
$$
This map is called the {\em R-trace}. As for ${\rm Tr}_{(i)}$, the
symbol $\Tr{i}$ stands for the R-trace over an $i$-th component in the tensor 
product
of the spaces $V$. Similarly, for the symbol $\Tr{i_1,i_2,\dots ,i_{L}}$.
\smallskip

Let us introduce special notation for the images of the elements 
$\kappa_i\in {\cal W}_n$ in the Yang--Baxter representation
$$ 
K_i:= \rho_R(\kappa_i)\, .
$$
By definition, the matrix $K=\rho_R(\kappa_1)\in  {\rm End}(V^{\otimes 2})$ is 
degenerate. As shown in \cite{IOP3}, rank of the operator $K$ corresponding to the 
skew invertible BMW type Yang--Baxter matrix  equals 1:
\be
\lb{rankK=1}
\mbox{rk}\,K=1.
\ee

The following set of relations for a skew invertible BMW type Yang--Baxter matrix  $R$ 
and 
for related matrices $D_R$, $K$, $\rho_R(a^{(i)})$ will be used in the main sections of 
the text. The proofs can be found in
\cite{OP-BMW}, sections 3.2, 3.4, and \cite{OP-SpCH}, section 2.4.
\ba
\lb{RDD}
R_{1}{D_R}_1 {D_R}_2 &=& {D_R}_1 {D_R}_2 R_{1}\, ,
\\[2pt]
\lb{KDD} 
K_{1} {D_R}_1 {D_R}_2 &=& {D_R}_1 {D_R}_2 K_{1}\ 
=\ \mu^2 K_{1}\, 
\\[2pt]
\lb{traceDK} 
\Tr{2} K_{1} &=& \mu\, I_1\, ,
\\[2pt]
\lb{traceD}
\tr_{\!\! R^{\phantom{A}\!\!\!}} I &=&
{\textstyle {(q-\mu)(q^{-1}+\mu)\over (q-q^{-1})}}\, ,
\\[2pt]
\lb{spec1}
\Tr{i} \rho_R(a^{(i)})& =& \delta_i(q,\mu)\, \rho_R(a^{(i-1)})\, ,
\\
\nn
 \mbox{where} \quad
\delta_i(q,\mu) &:=& -\, q^{i-1}\,{\textstyle{(\mu + q^{1-2i})(\mu^2 - q^{4-2i})\over
	(\mu + q^{3-2i})(q^i-q^{-i})  }}\, .
\ea

\subsection{Compatible Yang--Baxter pairs}
\lb{RF-pairs}

An ordered pair $\{ R, F\}$ of two Yang--Baxter matrices $R$ and $F$  is
called {\em a compatible Yang--Baxter pair} if the following conditions
\be 
R_1\, F_2\, F_1\, =\, F_2\, F_1\, R_2\, ,\qquad R_2\, F_1\, F_2\, =\, F_1\, F_2\, R_1\, ,
\lb{sovm}
\ee
are satisfied. The equalities (\ref{sovm}) are sometimes called {\em twist relations}. 
Clearly, $\{ R,P\}$ and $\{ R,R\}$ are compatible Yang--Baxter pairs.

\vskip .1cm
A compatible Yang--Baxter pair $\{R,F\}$ gives rise to a new Yang--Baxter matrix 
\be
\lb{R_f}
R_{_{\! F}}   := F^{-1} R F\, ,
\ee
called the {\em twisted} Yang--Baxter matrix. It is known that the Yang--Baxter pair 
$\{ R_{_{\! F}}   , F\}$ is again compatible. If $R$ is skew invertible and $F$ is strict 
skew 
invertible, then $R_{_{\! F}} $ 
is skew invertible; if additionally $R$ is strict skew invertible, 
then $R_{_{\! F}}  $ is strict skew invertible as well (for proofs, see \cite{OP-BMW}, proposition 
3.6). For a compatible pair $\{R,F\}$ of the skew invertible matrices the relation
\be
\lb{FCC} 
F_{1}\, {D_{R}}_1 {D_R}_2 \, =\, {D_R}_1 {D_R}_2 F_{1}\ ,
\ee
is satisfied (see \cite{OP-BMW}, corollary 3.4).

\section{Finite height and orthogonal Yang--Baxter matrices}
\lb{sec3}

\subsection{Finite height}
 In \cite{OP-SpCH} we have introduced a notion of the {\em finite height} BMW type 
 Yang--Baxter matrix. Namely, assuming that the parameters $q$ and $\mu$ of the 
 BMW type Yang--Baxter matrix  satisfy conditions (the union of 
 conditions (\ref{init-cond}), (\ref{q-mu-restrict}))
\be
\lb{q-mu-restrict-2} 
q\neq 0,\;\;q^{2i}\neq 1;  
\quad \mu\notin \{0, q\},\; \mu\neq -q^{-2i+1}\,  \;\;\forall\; i=1,\dots ,k\, ,
\ee
we consider Yang--Baxter matrix  realizations $\rho_R(a^{(i)})$ of the 
$q$-antisymmetrizers  (\ref{a^k}).
The skew invertible BMW type Yang--Baxter matrix  $R$ is called the 
{\em Yang--Baxter matrix  of  height $k$} if its $q$-antisymmetrizers satisfy conditions
\be
\lb{spec4} 
{\rm rk}\,\rho_R(a^{(i)})\neq 0\, \quad \forall\; i=2,3,\dots ,k, \quad\mbox{and} \quad
\rho_R\left( a^{(k)\uparrow 1}\sigma^-_1(q^{-2k})a^{(k)\uparrow 1} \right)\equiv 0\,
\ee
for some $k\geq 2$.

\smallskip

As discussed in \cite{OP-SpCH}, section 2.4, it is impossible to construct the height 
$k$ BMW type Yang--Baxter matrix  unless $\mu\in\{-q^{-1-2k},\pm q^{1-k}\}$.
Taking $\mu= -q^{-1-2k}$ one obtains series of the so called symplectic {\em Sp(2k) 
type} Yang--Baxter matrices. The corresponding quantum matrix algebras were 
investigated in \cite{OP-SpCH}. 

The choice $\mu= -q^{1-k}$ for even values of $k$ contradicts conditions 
(\ref{q-mu-restrict-2}). Besides, for odd values of $k$ it can be transformed to the case 
$\mu= q^{1-k}$ by the map: $R \mapsto  - R$. Therefore, a consideration of the case
$\mu= q^{1-k}$ completes the investigation of the finite height BMW type Yang--Baxter 
matrices and their corresponding QM-algebras. In the present work we are dealing 
with a family of the Yang--Baxter matrices belonging to the case $\mu= q^{1-k}$. 

\subsection{Orthogonal case}
\begin{defin}\lb{definition3.11}
Let $R$ be the skew invertible BMW type Yang--Baxter matrix  of the height $k\geq 2$ 
whose two
spectral values $q$ and $\mu$ are related by equality $\mu= q^{1-k}$.
Assume additionally that  
\be
\lb{rankA=1}
\mbox{rk}\,\rho_R(a^{(k)})=1.
\ee
Then $R$ is called an orthogonal Yang--Baxter matrix  of the type $O(k)$.
\end{defin}

\begin{rem}\lb{remark1.2}{\rm
We remind that for the $GL(k)$ type Yang--Baxter matrices, that is, the Hecke type 
Yang--Baxter matrices of the height $k$,   the strict skew invertibility implies  the 
relation (\ref{rankA=1})
\cite{Gu}. However, for the BMW type finite height Yang--Baxter matrices this is not the 
case.
E.g., for the family of symplectic $Sp(2k)$ type Yang--Baxter matrices, where the  
parameters of the BMW algebra are related by $\mu=-q^{-1-2k}$ (see \cite{OP-SpCH}, 
section 2.4),  the rank one condition is not valid for $k>1$. In the case $\mu=q^{1-k}$ a 
widely known family of Yang--Baxter matrices satisfying condition (\ref{rankA=1}) is 
related to the orthogonal quantum groups \cite{FRT}.  We believe that in this case the 
finiteness of the height together with the skew invertibility implies the rank one 
property. This is an interesting open question.}\hfill$\blacktriangleright$
\end{rem}
 
\subsection{Rank-one projectors}\lb{subsec3.5}

In this subsection we investigate specific
properties of certain rank-one projectors which can be constructed from the orthogonal 
Yang--Baxter matrices. Later on, these projectors will play a central role in the proof of 
the reciprocal relations (theorem \ref{theorem4.9}) and in the investigation of the 
quantum determinant (proposition \ref{proposition4.16a}).
\smallskip

First, we observe that any skew invertible BMW type Yang--Baxter matrix   gives rise 
to a series of
rank-one projectors.
\begin{lem}\lb{lemma3.13}
Assume that $R$ is a skew invertible Yang--Baxter matrix  of the BMW type. Then  
$\rho_R(c^{(2i)})\in {\rm End}(V^{\otimes 2i})$
(for the definition of the idempotents $c^{(2i)}$, see eqs.(\ref{kappa-i})$\,$) are 
rank-one operators fulfilling
relations
\ba\lb{trace-c2i}\Tr{2i}\rho_R(c^{(2i)}) &=& \eta^{-1}\mu\,  \rho_R(c^{(2i-2)\uparrow 1})
\quad \forall\; i\geq 1\, ,\\[1em]\lb{traces-c2i}\Tr{i+1,i+2,\dots ,2i}\rho_R(c^{(2i)}) &=&
(\eta^{-1}\mu)^i\, I^{\otimes i}\, .\ea\end{lem}

\nin{\bf Proof.~} 
For $i=1$, the relation (\ref{trace-c2i}) is just the equality (\ref{traceDK}). 
In case $i>1$, we check the equation (\ref{trace-c2i}) by the direct calculation
\ba\begin{array}{ccl} 
\Tr{2i} \rho_R(c^{(2i)}) &=& 
\eta^{-1}\,\rho_R(c^{(2i-2)\uparrow 1})
\left(\Tr{2i} K_{2i-1}\right) \rho_R\bigl( (\kappa_{2i-2}\kappa_{2i-3} \dots\kappa_{i+1})
(\kappa_1\kappa_2\dots \kappa_i)\bigr)\, 
\\[1em] 
&=&
\eta^{-1}\mu\, \rho_R\bigl( c^{(2i-2)\uparrow 1}
\kappa_{2i-2}\kappa_{2i-3}\dots\kappa_{i+1}(\kappa_1\kappa_2\dots \kappa_i)\bigr)\, 
\\[1em] 
&=&
\eta^{-1}\mu\, \rho_R\bigl( c^{(2i-2)\uparrow 1}\kappa_{i-1}
\dots\kappa_3\kappa_{2}(\kappa_1\kappa_2\dots
\kappa_i)\bigr)\, 
\\[1em] 
&=&
\eta^{-1}\mu\,  \rho_R\bigl( c^{(2i-2)\uparrow 1}\kappa_{i-1}\kappa_i
\bigr)\, =\, \eta^{-1}\mu\,  \rho_R\bigl(c^{(2i-2)\uparrow 1}\bigr)\, .
\end{array}
\nonumber
\ea
Here in the first line we substituted the expression (\ref{kappa-i2}) for $2i$-th order 
contractors. In the 
second line we used the relation (\ref{traceDK}). In the third line we applied repeatedly 
the
relations (\ref{idemp-c2}) for the elements $\kappa_{2i-2}, \dots , \kappa_{i+1}$ (they 
are quadratic polynomials in  $\sigma$'s); after each application of the 
relation (\ref{idemp-c2}) we move the resulting $\kappa$ to the right. In the last line we 
used
the relation (\ref{bmw5a}) $(i-2)$ times. Then, noticing that the element $\kappa_i$ 
divides
$c^{(2i-2)\uparrow 1}$ from the right (see (\ref{kappa-i2})) and applying once again eq.
(\ref{bmw5a}) we complete the transformation.

\smallskip
The relation (\ref{traces-c2i}) follows by a repeated application of the relation
(\ref{trace-c2i}).

\smallskip
As  we mentioned before, the rank one property of the 2-contractor 
$\rho_R(c^{(2)})=\eta^{-1} K$ associated 
to a skew invertible BMW type Yang--Baxter matrix  was proved in \cite{IOP3}.
To calculate the rank of the operators $\rho_R(c^{(2i)})\in {\rm End}(V^{\otimes 2i})$, 
$i>1$, we notice that using iteratively  
the relations (\ref{idemp-c1}) one can rewrite the formula (\ref{kappa-i}) for the 
contractors in the form
\be
\lb{kappa-i3} 
c^{(2i)}\, =\, \eta^{-\varepsilon(i)}\, c^{(2i-2)\uparrow 1} {\textstyle \left(\prod_{j=1}^{i} 
\kappa_{2j-1}\right)}
c^{(2i-2)\uparrow 1}\, ,
\ee
where $\varepsilon(i)$ is the parity function: $\varepsilon(i)=0/1$  if $i$ is even/odd. 
The appearance of the factor $\eta^{-1}$ for odd values of $i$ is due to the fact that 
for odd $i$ the iterative procedure terminates with 
$c^{(2)\uparrow (i-1)}=\eta^{-1}K_{i
}$. 
We leave further details to the reader.

By the rank-one property of the operator $K\in{\rm End}(V^{\otimes 2})$, a composite 
operator
$\rho_R(\prod_{j=1}^{i} \kappa_{2j-1}) = \prod_{j=1}^{i} K_{2j-1}\in 
{\rm End}(V^{\otimes 2i})$ also has rank one. Hence, the
operator $\rho_R(c^{(2i)})$ is either of rank-one or identically vanishes.
The latter possibility contradicts the relation (\ref{traces-c2i}). \hfill$\blacksquare$

\medskip
We next describe specific properties of the rank-one projectors.

\begin{lem}\lb{lemma3.14}
1) Let $R$ be a skew invertible Yang--Baxter matrix  of the BMW type.\footnote{Recall 
that 
for the BMW type the skew invertibility is equivalent to the strict skew invertibility.} 
Then
\ba\lb{spec-c1}
&& {\textstyle \left(\prod_{j=1}^{2i}{D_R}_j\right)}
\rho_R(c^{(2i)})\, =\, \mu^{2i}\, \rho_R(c^{(2i)})\quad \forall\, i\geq 1 .
\ea
Let $R$ be a Yang--Baxter matrix  of the type $O(k)$. Then
\ba\lb{spec-a1}
&&
{\textstyle \left(\prod_{j=1}^{k}{D_R}_j\right)} \rho_R(a^{(k)})\, =\, q^{k(1-k)}\, 
\rho_R(a^{(k)})\, .
\ea

\smallskip
2) Let additionally $F$ be a strict skew invertible Yang--Baxter matrix  and assume that 
$\{R,F\}$
is a compatible pair. Then
one has, respectively,
\ba
\lb{spec-c2} 
&& {\textstyle \left(\prod_{j=1}^{2i}{D_{(R_{_{\! F}}  )}}_j\right)}
\rho_R(c^{(2i)})\, =\, \mu^{2i}\, \rho_R(c^{(2i)})\quad \forall\, i\geq 1
\ea
if the Yang--Baxter matrix  $R$ is strict skew invertible of the BMW type, and
\ba\lb{spec-a2} &&
{\textstyle\left(\prod_{j=1}^{k}{D_{(R_{_{\! F}}  )}}_j\right)}\rho_R(a^{(k)})\, =\, q^{k(1-k)}\, 
\rho_R(a^{(k)})\, \ea
if $R$ is of the type $O(k)$.
\end{lem}

\nin {\bf Proof.~} We first prove the equality (\ref{spec-c1}). Substitute the expression 
(\ref{kappa-i3}) into the left hand side of (\ref{spec-c1}).
By the left equality in (\ref{KDD}), we can change the order of operators
$\rho_R(c^{(2i-2)\uparrow 1})$ and $\prod_{j=1}^{2i}{D_R}_j$. 
Then, by the right equality in (\ref{KDD}), the operator 
$\prod_{j=1}^{2i}{D_R}_j$
in the product 
$\prod_{j=1}^{2i}{D_R}_j\,\cdot\,\rho_R(\prod_{j=1}^{i} \kappa_{2j-1})$ 
can be replaced by the factor $\mu^{2i}$. 

\medskip
The equality (\ref{spec-a1}) is confirmed by the calculation
\begin{eqnarray}
\nonumber
\textstyle \left(\prod_{j=1}^{k}{D_R}_j\!\right)\rho_R(a^{(k)})\;\,=
& \!\!
\rho_R(a^{(k)})
\left(\prod_{j=1}^{k}{D_R}_j\!\right) 
\rho_R(a^{(k)})
&\!\!		=\,
\left(\Tr{1,2,\dots k}\rho_R(a^{(k)})\right)\rho_R(a^{(k)}) 
\\[2pt]
\nonumber
&\, =\, q^{k(1-k)}\, \rho_R(a^{(k)}). &
\end{eqnarray}
Here the idempotency of the operator $\rho_R(a^{(k)})$ and its commutativity with the 
product
$\prod_{j=1}^k {D_R}_j$ (by the relation
(\ref{RDD})) were taken into account in the first equality; the second equality follows
by the rank-one property of the operator $\rho_R(a^{(k)})$. Eqs.(\ref{spec1})
were used for the evaluation of the R-traces in the second line. Namely,
\be
\lb{qdim-O(k)}
\Tr{1,2,\dots ,i}\rho_R(a^{(i)}) = \prod_{j=1}^i \delta_j(q,\mu)|_{\mu=q^{1-k}}\, =\,
q^{-ki}\,\, {q^{2i}+q^{k}\over 1+q^k}\, {k_q!\over i_q! (k-i)_q!},
\ee
where $j_q!:=1_q 2_q 3_q\dots j_q$. 

\smallskip
The conditions of the second part of the proposition guarantee that the twisted 
Yang--Baxter matrix  $R_{_{\! F}}  =F^{-1}RF$ is 
skew invertible of the BMW type. 
Hence, the relations (\ref{spec-c1}) and (\ref{spec-a1}) are 
satisfied for $R_{_{\! F}}  $ as well.
Now, to prove the equalities (\ref{spec-c2}) and (\ref{spec-a2}), it is enough to 
demonstrate that the
two operators $\rho_{R_{_{\! F}} }(c^{(2i)})$ and $\rho_R(c^{(2i)})$ (respectively, 
$\rho_{R_{_{\! F}}  }(a^{(k)})$
and $\rho_R(a^{(k)})$) are related by a similarity transformation, which commutes with
$\prod_{j=1}^{2i} {D_{(R_{_{\! F}}  )}}_j$ (respectively, $\prod_{j=1}^k {D_{(R_{_{\! F}}  )}}_j$).

\smallskip
Define a set of operators $Z^{(i)}\in {\rm Aut}(V^{\otimes i})$, $i=1,2,\dots ,$ by the 
iteration\footnote{Two expressions for the element $Z^{(i+1)}$ in eq.(\ref{Z-i}) 
correspond to two different decompositions of the lift, to the braid group 
${\cal B}_{i+1}$, of the longest element of the symmetric group $S_{i+1}$.}
\be
\lb{Z-i} 
Z^{(1)}:=I, \qquad Z^{(2)}:=F_1,\qquad Z^{(i+1)}:=\left(F_1 F_2 \dots F_i\right) Z^{(i)} =
Z^{(i)} \left(F_i \dots F_2 F_1\right) .
\ee

By induction, one can readily prove the relations
\be\lb{Z-R} R_i\, Z^{(k)}\, =\, Z^{(k)}\, {(R_{_{\! F}}  )}_{(k-i)}\,
\qquad \forall\; i=1,2,\dots ,k-1\, .
\ee

Namely, for $i=2,\dots ,k-1$, one has, by induction on $k$,
\ba\nonumber R_i Z^{(k)} &=& R_i F_1 F_2 \dots F_{k-1} Z^{(k-1)} = F_1 F_2 
\dots F_{k-1} R_{i-1}Z^{(k-1)}
\\[1em]\nonumber &=&F_1 F_2 \dots F_{k-1} Z^{(k-1)}{(R_{_{\! F}} )}_{(k-i)} = 
Z^{(k)}{(R_{_{\! F}} )}_{(k-i)}\ .\ea
Here we used the twist relation for the compatible pair $\{ R,F\}$ and the first iterative
definition of the element $Z^{(k)}$ in eq.(\ref{Z-i}).

{}For $i=1$, one has, again by induction on $k$,
\ba\nonumber R_1 Z^{(k)} &=& R_1 Z^{(k-1)} F_{k-1}\dots F_2 F_1 = 
Z^{(k-1)}{(R_{_{\! F}} )}_{(k-2)}F_{k-1}\dots F_2 F_1
\\[1em]\nonumber &=&Z^{(k-1)}F_{k-1}\dots F_2 F_1{(R_{_{\! F}} )}_{(k-1)}=
Z^{(k)}{(R_{_{\! F}} )}_{(k-1)}\ .\ea
Here we used the twist relation for the pair $\{ R_{_{\! F}} ,F\}$ and the second iterative 
definition
of the element $Z^{(k)}$ in eq.(\ref{Z-i}).

\medskip
The relations (\ref{Z-R}) imply following equalities
\be
\lb{Z-ac1} 
\rho_R(a^{(k)})\, Z^{(k)} =  Z^{(k)} \rho_{R_{_{\! F}} }(a^{(k)})\, ,
\ee
\be
\lb{Z-ac2} 
\rho_R(c^{(2i)})\,  Z^{(2i)} =  Z^{(2i)} \rho_{R_{_{\! F}} }(c^{(2i)})\, 
\quad\forall\; i=1,2,\dots
\ee
To verify eq.(\ref{Z-ac1}), substitute different expressions from eq.(\ref{a^k})
for the $q$-antisymmetrizers, one expression for $\rho_R(a^{(k)})$ and 
another one for $\rho_{R_{_{\! F}} }(a^{(k)})$; then use
eq.(\ref{Z-R}). As for eq.(\ref{Z-ac2}), note that the expression (\ref{kappa-i}) 
for the contractor $c^{(2i)}$ is invariant with
respect to the reflection $\sigma_j\leftrightarrow\sigma_{2i-j}$. Thus
it can be used in both sides of eq.(\ref{Z-ac2}).

\smallskip
Due to the relations (\ref{FCC}) for the compatible pair
 $\{R_{_{\! F}},F\}$, the operators $Z^{(l)}$ and 
 $\prod_{j=1}^l {D_{(R_{_{\! F}} )}}_j$
commute for any $l$ 
and, hence, the operator $Z^{(2i)}$ (respectively, the operator $Z^{(k)}$) 
realizes the
similarity transformation we are looking for.\hfill$\blacksquare$

\begin{rem}{\rm For a strict skew invertible Yang--Baxter matrix  of the Hecke 
type of height $k$ the identity (\ref{spec-a1}) holds as well. If the Yang--Baxter matrix  
$R$ is strict skew invertible of the Hecke type, 
$F$ is strict skew invertible Yang--Baxter matrix  and  $\{R,F\}$
is a compatible pair then the identity (\ref{spec-a2}) also holds. 
The above proof holds in the Hecke case. } \hfill$\blacktriangleright$
\end{rem}
\vskip .8cm

\section{Quantum matrix algebra in the orthogonal case}\lb{sec4}

In this section we introduce the main objects of our study, the {\em quantum matrix 
(shortly, QM-) algebras} of the orthogonal type.
\smallskip 

In subsections \ref{subsec4.1}-\ref{subsec4.1b} we  collect some known results about the BMW type 
QM-algebras, in particular, about the structure of their distinguished 
commutative subalgebra, 
called the {\em characteristic subalgebra}. Additional {\em reciprocal} relations,
specific to the orthogonal type QM-algebras,  
satisfied by generators of the characteristic subalgebra  are 
derived in subsection \ref{subsec4.11} and resolved in subsection \ref{components}. 
Subsections \ref{qdetorth} and \ref{repqma} present results, respectively, on the 
properties of the quantum determinant (analogue of the determinant for the quantum 
matrix) and on two specific representations of the orthogonal QM-algebra. These results 
are used in the solution of the reciprocal relations.

\subsection{Orthogonal quantum matrix algebra}\lb{subsec4.1}

Here we remind necessary definitions and few basic results about quantum matrix 
algebras of BMW type. For a detailed exposition the reader is referred to 
\cite{OP-BMW}. We remind that $V$ is a finite dimensional ${\Bbb C}$-linear space, 
$\dim V = \mbox{\sc n}$;  the operators $R$ and $F$ belong to 
${\rm Aut}(V^{\otimes 2})$.

\begin{defin}\lb{definitionQMA}
Assume that $\{R,F\}$ is a compatible pair of the strict skew invertible 
Yang--Baxter 
matrices.
The QM-algebra ${\cal M}(R,F)$ is a quotient algebra of the free associative unital 
algebra $W={\Bbb C}\langle M^a_b\rangle_{a,b=1}^{\mbox{\footnotesize\sc n}}$ 
by a two-sided ideal generated by 
entries of the matrix relation
\be 
R_1 M_{\overline 1}M_{\overline 2} = M_{\overline 1}M_{\overline 2}R_1\ . 
\label{qma}
\ee
Here $M = \|M^a_b\|_{a,b=1}^{\mbox{\footnotesize\sc n}}$ is the matrix of 
generators;
the matrix copies $M_{\overline i}$, $i\geq 1$, are constructed by the following 
iteration 
\be
M_{\overline 1}:=M_1, \quad M_{\overline{i}}:=
F^{\phantom{-1}}_{i-1}M_{\overline{i-1}}F_{i-1}^{-1}\ .
\lb{kopii}
\ee	
	
The QM-algebra ${\cal M}(R,F)$ is
called BMW$/O(k)$ type if, respectively,   
the Yang--Baxter matrix  $R$ is of BMW$/O(k)$ type. 	
\end{defin}

We call the matrix $M$ whose entries are the generators of the QM-algebra 
${\cal M}(R,F)$ the {\em quantum matrix}.
\smallskip

For any given value of the index $i\geq 1$ the set of relations 
\be 
R_i M_{\overline i}M_{\overline{i+1}} = M_{\overline i}M_{\overline{i+1}}R_i 
\label{qmai}
\ee
is equivalent to (\ref{qma}) (see \cite{OP-BMW}, lemma 4.3).\smallskip

For the BMW type QM-algebra the relations 
\be
\lb{tau2}
K_{i}\, M_{\overline{i}}M_{\overline{i+1}} =
M_{\overline{i}}M_{\overline{i+1}}\, K_{i}\,\, =\,\,\mu^{-2}g\,  K_{i} 
\quad 
\forall\; i\geq 1 , 
\ee
where
\be\lb{tau}
g := \mu^2\,	\tr_{(1,2)} \left( M_{\overline{1}}M_{\overline{2}}\,  \rho_R(c^{(2)})\right) 
\, =\, \eta^{-1}\,\Tr{1,2} \Bigl( M_{\overline{1}}M_{\overline{2}}\, K_1\Bigr) ,
\ee
are satisfied as a consequence of eqs. (\ref{qmai}), (\ref{KDD}) and the rank-one 
property of the matrix $K$ (see \cite{OP-BMW}, lemma 4.6).
The element $g$ is called the  {\em contraction of $M$}. 

\begin{rem} {\rm
(i) For orthogonal group $O(n)$, $n\geq 2$, the tensor square of the vector representation $V$
decomposes into the direct sum of three irreducible representations for each $n$ except 
$n=4$. For $n=4$ the tensor square of the vector representation is the direct sum of four 
irreducible representations: the wedge square of $V$ is a direct sum of self-dual 
and anti-self-dual tensors. This geometry quantizes and there are four $q$-deformed projectors
(see \cite{OSWZ1, OSWZ2}). This is also explained by the $q$-deformation of 
accidental isomorphisms of semi-simple Lie groups \cite{JO}. In such situation the QM-algebra 
must be defined by the union of equations 
\be\lb{genqma} \Pi_\alpha M_{\overline 1}M_{\overline 2} = 
M_{\overline 1}M_{\overline 2} \Pi_\alpha\ee 
for all projectors $\Pi_\alpha$ $q$-deforming the decomposition of the tensor square 
of $V$ into the direct sum of irreducibles. 

\vskip .4cm
(ii) We believe that for any $O(4)$-type Yang--Baxter matrix, admitting the classical limit, 
there exists a non-trivial
decomposition of the second $q$-antisymmetrizer (acting in the tensor square of the vector 
representation) into a sum of two projectors, generalizing the 
situation with the standard $q$-deformation of $O_4$. 

\vskip .4cm
(iii) The fact that in general the QM-algebra is a quotient of the algebra defined by the eq.(\ref{qma})
does not affect the main objective of our study -- the characteristic equation 
for the quantum matrices. 
This should 
be compared with the perfectly classical situation -- the characteristic polynomial $\chi_A(t)$ 
of a matrix $A$ is always just $\chi_A(t):=\rdet (tI-A)$, regardless of the type of the matrix 
$A$ (general linear, symplectic, orthogonal, ...) Note however that the type of the matrix 
affects the properties of its eigenvalues.

\vskip .4cm
(iv) Let $g$ be a semi-simple Lie algebra, $Ug$ its universal enveloping algebra and $U_q g$ 
its standard $q$-deformation. Let $\rho$ be an irreducible representation of $g$ in a vector space $V$,
and $\rho_q$ 
its $q$-deformation. The image of the universal $\cal{R}$-matrix may not distinguish 
all the projectors 
in the tensor square of $V$, as it is the case for $O(4)$. In such a situation the quantum matrix 
algebra should be defined by the union of equations (\ref{genqma}) as above. However, in the 
$O(4)$ case the situation is quite beautiful: one can construct, out of projectors $\Pi_\alpha$ 
two Yang--Baxter matrices \cite{OSWZ1} 
(which, of course commute, since they are built of the projectors 
entering the decomposition of unity). These two matrices fully serve the geometry of the 
$q$-deformed 4-dimensional orthogonal space. We believe that this will be the case for all 
$g$ and $\rho$ for which the image of the universal $\cal{R}$-matrix does not distinguish 
all the projectors 
in the tensor square of $V$ (maybe there will be more than two Yang--Baxter matrices) but this is an 
open question. \hfill$\blacktriangleright$}
 \end{rem}

\subsection{Characteristic subalgebra}\lb{subsec4.1b}

\begin{def-prop} {\em \cite{IOP1}}
Let ${\cal C}(R,F)$ be a subspace in ${\cal M}(R,F)$ which is
spanned linearly by the unity and elements
\be
\lb{char}
{\rm ch}(\alpha^{(n)}) := \Tr{1,\dots ,n}(M_{\overline 1}\dots M_{\overline n}\,
\rho_R(\alpha^{(n)}))\ ,\quad n =1,2,\dots\ ,
\ee
where $\alpha^{(n)}$ is an arbitrary element
of the braid group ${\cal B}_n$. The space ${\cal C}(R,F)$ is a commutative 
subalgebra in ${\cal M}(R,F)$. It is called the  characteristic subalgebra of 
${\cal M}(R,F)$.
\end{def-prop}

\begin{rem}
\lb{rem3}
 {\rm
For the family  ${\cal M}(R,R)$ of the so-called {\em reflection equation} algebras  
their characteristic subalgebras lie in their centers: 
${\cal C}(R,R)\subset Z[{\cal M}(R,R)]$   (see \cite{I}, proposition 5).}
\hfill$\blacktriangleright$
\end{rem}

\begin{prop}\lb{gen-char}	{\em(\cite{OP-BMW}, propositions 4.7, 4.8)}
The characteristic subalebra of the BMW type QM-algebra is generated by either of 
the sets:
\begin{enumerate}
\item $\{g,\goth{p}_i\}_{i\geq 0}$, where elements
\be
\lb{power sums}
\goth{p}_0:=\tr_{\!\! R}\, I=\mu \eta, \quad \goth{p}_1 := \tr_{\!\! R}\, M\, ,\quad
\goth{p}_i :=  {\rm ch}(\sigma_{i-1}\dots\sigma_2\sigma_1) 
,  \; \forall\,i\geq 2,
\ee
are called {\it power sums};
\item  $\{g,\goth{h}_i\}_{i\geq 0}$, where elements
\be
\lb{complete symm-f}
{\goth{h}_0} :=1, \quad
{\goth{h}_i} := {\rm ch}(s^{(i)}), \; \forall\, i\geq 1,
\ee
are called {\it complete  sums};

\item $\{g,\goth{e}_i\}_{i\geq 0}$, where elements
\be
\lb{elementary symm-f}
\goth{e}_0 :=1, \quad
\goth{e}_i := {\rm ch}(a^{(i)}), \; \forall\, i\geq 1,
\ee
are called {\it elementary sums}.
\end{enumerate}
Here in the definitions of the last two sets we assume additionally that  conditions 
(\ref{q-mu-restrict}) on the parameters $q$ and $\mu$ are fulfilled.
\end{prop}

Our definition of the QM-algebra does not assume invertibility of the 
quantum matrix $M$ with respect to the usual matrix product. 
Now we introduce a natural extension of the BMW type QM-algebra 
containing the inverse 
of the matrix $M$.

\begin{def-prop}
\lb{proposition4.11} 
{\em(\cite{OP-BMW}, lemma 4.13 and proposition 4.14)} Let $
{\cal M^{^\bullet\!}}(R,F)$ and ${\cal C^{^\bullet\!}}(R,F)$ denote  
localizations of the BMW type QM-algebra 
${\cal M}(R,F)$ and of its characteristic 
subalgebra ${\cal C}(R,F)$
by adding the inverse  $g^{-1}$ of the contraction 
\be
\lb{j-inv} 
g^{-1}\, g\, =\, g\, g^{-1}\, =\, 1\, .
\ee
The contraction and its inverse satisfy the following permutation relations with the 
QM-algebra gene\-rators
\be
\lb{g-perm}
M\, g\, =\, g\, ( G^{-1}M G)\, ,\ \ \mbox{and so}\ \ \ 
 g^{-1}\, M\, =\, (G^{-1}MG)\, g^{-1}\, ,
\ee	
where the numeric mutually inverse matrices $G^{\pm 1}\in {\rm Aut}(V)$ are given by 
the formulas
\be
\lb{G}
G_1 \, :=\, \tr_{(23)} K_2^{} F_1^{-1} F_{2}^{-1}\, , \qquad  
G_1^{-1}\, =\, \tr_{(23)} F_2 F_1 K_2\, .
\ee
The  algebra ${\cal M^{^\bullet\!}}(R,F)$ contains
the inverse to the quantum matrix $M$,
\be
\lb{M-inv} 
(M^{-1})_1\, =\, \mu\, \Tr{2}\! \left( 
M_{\overline 2}  K_{1}^{}\right)\, g^{-1},
\ee 
so that
$$
M\, M^{-1}\, =\, M^{-1} M\, =\, I.
$$
\end{def-prop}

\subsection{Reciprocal relations}\lb{subsec4.11}

For the QM-algebra ${\cal M}(R,F)$ where $R$ is of finite height, 
the set of the elementary sums becomes finite. 
Moreover, for the  orthogonal $O(k)$ type QM-algebra, the rank-one condition  
(\ref{rankA=1}) for the operator
$\rho_R(a^{(k)})$
gives rise to algebraic equalities among the elementary sums
$\{ \goth{e}_i\}_{0\leq i\leq k}$ and the contraction $g$. These equalities --- we call them 
{\em reciprocal relations} --- are described in the next theorem.

\begin{theor}\lb{theorem4.9}
Let the QM-algebra ${\cal M}(R,F)$ be of the orthogonal $O(k)$ type. Then
the following reciprocal relations for the generators $\{g,\goth{e}_i\}_{i=0}^k$ 
of its characteristic subalgebra
\be
\lb{reciprocal} 
g^{k-i}\, \goth{e}_i\, =\, \goth{e}_k\, \goth{e}_{k-i}\,\qquad \forall\; i=0,1,\dots ,k-1, 
\ee
are satisfied.
\end{theor}

\nin {\bf Proof.~} For the proof, we use three auxiliary statements.

\begin{lem}\lb{lemma4.10}
Let the QM-algebra ${\cal M}(R,F)$ be of the BMW type, respectively, of the $O(k)$ 
type. 
Then we have, respectively,
\ba\lb{spec-cM} && M_{\overline{1}}\, M_{\overline{2}}\dots M_{\overline{2i}}\,\,
\rho_R(c^{(2i)})\, =\, (\mu^{-2} g)^i\, \rho_R(c^{(2i)})\, \ea
in the BMW case and
\ba\lb{spec-aM} && M_{\overline{1}}\,M_{\overline{2}}\dots M_{\overline{k}}\,\,
\rho_R(a^{(k)})\, =\, q^{k(k-1)}\, \goth{e}_k\, \rho_R(a^{(k)})\, \ea
in the $O(k)$ case. 
\end{lem}

\nin {\bf Proof.~} The relations (\ref{spec-cM}) and (\ref{spec-aM}) are proved in 
exactly the
same way as the relations (\ref{spec-c1}) and (\ref{spec-a1}) of the lemma 
\ref{lemma3.14}.
One has to use additionally the relation (\ref{tau2}) for the proof of the relation
(\ref{spec-cM}).
{}For the relation (\ref{spec-aM}), the equality
\be
\lb{qdet-1}
\phantom{a}\hspace{-5mm}
\tr_{(1,2,\dots k)} (M_{\overline{1}}\,M_{\overline{2}}\dots M_{\overline{k}}
\,\rho_R(a^{(k)}))\, =\, q^{k(k-1)}\,\Tr{1,2,\dots k}
(M_{\overline{1}}\,M_{\overline{2}}\dots M_{\overline{k}} \,
\rho_R(a^{(k)}))\, =\, q^{k(k-1)}\, \goth{e}_k ,
\ee
following from the relation (\ref{spec-a1}), also should be taken into account.
\hfill$\blacktriangleright$

\begin{rem}{\rm For a strict skew invertible Yang--Baxter matrix  of the Hecke 
type of height $k$ the identity (\ref{qdet-1}) holds as well, with 
the same proof. } \hfill$\blacktriangleright$
\end{rem}

\begin{lem}\lb{lemma4.11}
Let $Z^{(i)}\in {\rm Aut}(V^{\otimes i})$ be the operators defined in the formulas 
(\ref{Z-i}). For a given matrix
$M\in {\rm Mat}_{\mbox{\footnotesize\sc n}}(W)$, consider yet another set of 
its "copies" defined as (c.f.
with eq.(\ref{kopii})$\,$)
\be\lb{kopii-2} M_{\underline{1}}\, :=\, M_1 ,\qquad M_{\underline{i+1}}\, :=\,
{}F_i^{-1}M_{\underline{i}}\, F_i , \qquad i=2,3,\dots .\ee
Then relations
\be \lb{Z-M}M_{\overline{j}}\, Z^{(i)}\, =\, Z^{(i)}\, M_{\underline{i-j+1}}\ee
hold for any $i=1,2,\dots$ and for any $j=1,2,\dots ,i$.\end{lem}

\nin {\bf Proof.~} The equalities (\ref{Z-M}) can be proved by induction. First, one 
considers the case
$j=i$ and carries out induction on $i$. Then, performing induction on $(i-j)$, one 
proves
the equalities (\ref{Z-M}) in the full generality. \phantom{a}\hfill$\blacktriangleright$

\begin{lem}\lb{ajc2j} The following relation
\be a^{(j)\uparrow j}\, c^{(2j)}\, =\, a^{(j)}\, c^{(2j)}\lb{noca1}\ee
holds.
\end{lem}

\nin {\bf Proof.~}
The proof is by induction on $j$ (the base of induction, $j=1$, is obvious):
\be\begin{array}{l}  a^{(j+1)\uparrow j+1}c^{(2j+2)}
= {\displaystyle\frac{q^j}{(j+1)_q}}a^{(j)\uparrow j+1}\sigma_{2j+1}^-(q^{-2j})\,
a^{(j)\uparrow j+1}c^{(2j+2)} \\[1em]
\ \ \ \ \
= {\displaystyle\frac{q^j}{(j+1)_q}}a^{(j)\uparrow j+1}\sigma_{2j+1}^-(q^{-2j})
a^{(j)\uparrow 1}c^{(2j+2)}
= {\displaystyle\frac{q^j}{(j+1)_q}}a^{(j)\uparrow j+1}a^{(j)\uparrow 1}
\sigma_{1}^-(q^{-2j}) c^{(2j+2)}\\[1em]
\ \ \ \ \
={\displaystyle\frac{q^j}{(j+1)_q}}a^{(j)\uparrow 1}
\sigma_{1}^-(q^{-2j}) a^{(j)\uparrow 1}c^{(2j+2)}=a^{(j+1)}c^{(2j+2)}\ .
\end{array}\lb{noca2}\ee
In the first equality we used the left definition of the 
$q$-antisymmetrizers in eq.(\ref{a^k});
in the second equality we used the induction assumption (the contractor 
$c^{(2j)\uparrow 1}$
divides the contractor $c^{(2j+2)}$, see eq.(\ref{idemp-c1}));
in the third equality we moved the element $\sigma_{2j+1}^-(q^{-2j})$ rightwards to the
contractor and used the relation (\ref{idemp-c2}); in the fourth equality we moved
the $q$-antisymmetrizer $a^{(j)\uparrow j+1}$ rightwards to the contractor and again
used the induction assumption; in the fifth equality we used the right definition of the
$q$-antisymmetrizers in eq.(\ref{a^k}).\hfill$\blacktriangleright$

\medskip
Now, we will prove the relations (\ref{reciprocal}) by expressing elements 
${\rm ch}(c^{(2k-2i)}a^{(k)\uparrow k-i})
\in {\cal C}(R,F)$ in terms of elementary generators in two different ways.

\smallskip
The first way is
\be\lb{finoca}
\begin{array}{l}
{\rm ch}( c^{(2k-2i)} a^{(k)\uparrow k-i})
\\[1em]
\ \ \ \ \ = (\mu^{-2} g)^{k-i}\, \Tr{1,2,\dots ,2k-i}
\left( M_{\overline{2k-2i+1}}\dots M_{\overline{2k-i-1}} M_{\overline{2k-i}}\,\rho_R
(c^{(2k-2i)}a^{(k)\uparrow k-i})\right)
\\[1em]
\ \ \ \ \  =(\mu^{-2} g)^{k-i}\,\Tr{1,2,\dots ,2k-i}
\left(Z^{(2k-i)}M_{\underline{i}}\dots M_{\underline{2}} M_{\underline{1}}\,\,
\rho_{R_{_{\! F}} }(c^{(2k-2i)
\uparrow i} a^{(k)}) (Z^{(2k-i)})^{-1}\right)
\\[1em]
\ \ \ \ \ =(\mu^{-2} g)^{k-i}\,
\Tr{1,2,\dots ,i} M_{\underline{i}}\dots M_{\underline{2}} M_{\underline{1}}\,
\Bigl( \TR{i+1,\dots ,2k-i}
{R_{_{\! F}} }\, \rho_{R_{_{\! F}} }(c^{(2k-2i)\uparrow i} a^{(k)})\Bigr) \\[1em]
\ \ \ \ \ =(\eta^{-1}\mu^{-1} g)^{k-i}\,
\Tr{1,2,\dots ,i} M_{\underline{i}}\dots M_{\underline{2}} M_{\underline{1}}\, 
\Bigl( \TR{i+1,\dots ,k}
{R_{_{\! F}}  }\,\rho_{R_{_{\! F}}  }(a^{(k)})\Bigr) \\[1em]
\ \ \ \ \ ={\displaystyle \bigl(\prod_{j=1}^{k-i}\,\eta^{-1}\mu^{-1}
\delta_{k-j+1}\, g\bigr)}\,\Tr{1,2,\dots ,i}\left( (Z^{(i)})^{-1} M_{\overline{1}}
M_{\overline{2}}\dots
M_{\overline{i}}\, \rho_R(a^{(i)}) Z^{(i)}\right) \\[1.5em]
\ \ \ \ \ ={\displaystyle\bigl(\prod_{j=1}^{k-i}\,
\eta^{-1}\mu^{-1}\delta_{k-j+1}\bigr)}\, g^{k-i}\, \goth{e}_i\, .\end{array}\ee
Here the relation (\ref{spec-cM}) was applied  in the first equality; the
relations (\ref{Z-M}) as well as direct generalizations
\ba
\nonumber
\rho_R(c^{(2i)}) Z^{(j)}& =& Z^{(j)} \rho_{R_{_{\! F}}  }(c^{(2i)\uparrow j-2i})\qquad 
\forall\; j\geq 2i\geq 0\, ,
\\
\nonumber
\rho_R(a^{(i)\uparrow j-i}) Z^{(j)}& =&  Z^{(j)}\rho_{R_{_{\! F}}  }(a^{(i)})\qquad\qquad\; 
\forall\; j\geq i\geq 0
\ea
of the relations (\ref{Z-ac1}-\ref{Z-ac2}) were used 
in the second equality; the commutativity of 
$\prod_{j=1}^{2k-i}{D_R}_j$ with $Z^{(2k-i)}$ (which follows from the relations 
(\ref{FCC})\,)
allows to cancel the latter when passing to third equality; the relations (\ref{spec-c1}) 
and  (\ref{spec-c2}) were used
to substitute the R-trace $\tr_{\!\!\raisebox{-3pt}{$\scriptstyle R$}}$ by the R-trace
$\tr_{\!\!\raisebox{-1.5pt}{$\scriptstyle R_{_{\! F}}  $}}$ in the spaces with 
labels ~$\scriptstyle i+1,i+2,\dots ,
2k-i$~ in the third equality; the relations (\ref{traces-c2i})  and (\ref{spec1}) were 
applied for the evaluation of the R-traces
in the fourth and fifth equalities; a similarity transformation by $Z^{(i)}$ was performed 
in the argument of the
R-trace and, then, the relations (\ref{Z-M}) and (\ref{Z-ac1}-\ref{Z-ac2}) were used to 
complete the transformation.

\smallskip
The second way to express the element
${\rm ch}(c^{(2k-2i)}a^{(k)\uparrow k-i})$ in terms of elementary generators is
\be\begin{array}{l} {\rm ch}(\lefteqn{ c^{(2k-2i)} a^{(k)\uparrow k-i}) =
q^{k(k-1)}\, \Tr{1,2,\dots ,2k-i}\Bigl(M_{\overline{1}} M_{\overline{2}}\dots 
M_{\overline{k-i}}\,
\rho_R(a^{(k)\uparrow k-i} c^{(2k-2i)})\Bigr)\, \goth{e}_k } \\[1em] =q^{k(k-1)}\,
\bigl({\displaystyle \prod_{j=k-i+1}^{k} \delta_j}\bigr)\, \Tr{1,2,\dots ,2k-2i}
\Bigl( M_{\overline{1}}
M_{\overline{2}}\dots M_{\overline{k-i}}\,\rho_R(a^{(k-i)} c^{(2k-2i)})\Bigr)\, \goth{e}_k 
\hspace{10mm}\\[1.5em]
=q^{k(k-1)}\, \bigl({\displaystyle \prod_{j=k-i+1}^{k} \delta_j}\bigr)\,
\bigl( \eta^{-1}\mu\bigr)^{k-i}\ \Tr{1,2,\dots ,k-i}\Bigl( M_{\overline{1}}
M_{\overline{2}}\dots M_{\overline{k-i}}\,\rho_R(a^{(k-i)})\Bigr) \, \goth{e}_k \\[1em] =
\bigl({\displaystyle \prod_{j=1}^{k-i} \eta^{-1}\mu\, \delta_j^{-1}}\bigr)\,\goth{e}_{k-i}\, \goth{e}_{k}\, .
\end{array}\lb{noca}\ee
Here (a version of) the relation (\ref{spec-aM}) was used in the first line; eq.
(\ref{spec1})
was applied for the evaluation of the R-traces in the second line. There also, we
substituted, due to the lemma \ref{ajc2j}, the element 
$a^{(k-i)\uparrow k-i}\, c^{(2k-2i)}$
by the element $a^{(k-i)} c^{(2k-2i)}$ in the argument of $\rho_R$; in the third line,
some R-traces were evaluated with
the help of the relation (\ref{traces-c2i}); the relation  (\ref{qdim-O(k)}) for $i=k$ were 
used to transform the expression to its final form in the fourth line.

\medskip
Comparing the results of the calculations (\ref{finoca}) and (\ref{noca}) and substituting 
the
values, which $\mu$ and $\delta_j$ take in
the $O(k)$ case (see the definition \ref{definition3.11} and eq.(\ref{spec1})$\,$), we 
finish
the verification of the relations (\ref{reciprocal}).
\hfill$\blacksquare$

\subsection{Group-like elements in the orthogonal QM-algebras}\lb{qdetorth}

A particular subfamily of the QM-algebras --- ${\cal M}(R,P)$ --- is often called {\em FRT-algebras} as it was first considered in details in \cite{FRT}. These algebras carry a natural bialgebra structure with the coproduct given on generators by$\,$\footnote{In this article summation is always assumed to be performed from lower to upper indices; that is, for example, 
if $A$ and $B$ are two endomorphisms of a vector space $V$, then
$(AB)^i_j:=\sum_l A^i_l B^l_j$.}%
\be
\lb{coproduct}
\Delta (M^a_b) = \sum_c M^a_c\otimes M^c_b.
\ee 
In this subsection we investigate algebraic properties of the group like elements (that is, elements 
$b$ such that $\Delta  b = b\otimes b$) in the FRT-algebras. A generalization of the coproduct 
(\ref{coproduct}) for a general QM-algebra ${\cal M}(R,F)$ is discussed in 
the appendix \ref{appendix-A}.
\medskip

Let $b^{(i)}\in {\cal B}_i$ be such that $\rho_R(b^{(i)})\in{\rm End}\left(V^{\otimes i}\right)$ 
is a rank-one projector. Then the element $b\in {\cal M}(R,P)$, defined by 
\be
\lb{group-like}
b = \tr_{(1,2,\dots ,i)} (M_{\overline{1}}\,M_{\overline{2}}\dots M_{\overline{i}}
\,\rho_R(b^{(i)})) 
\ee
is group-like. 

\medskip 
Thus, in any BMW type FRT-algebra ${\cal M}(R,P)$ the rescaled contraction $\mu^{-2} g$ associated with the rank-one projector $\rho_R(c^{(2)})$ is group-like.\footnote{ The series of the higher order rank-one projectors  $\rho_R(c^{(2i)})$, $i>1$, do not give rise to the new group-like elements. The projector $\rho_R(c^{(2i)})$ is associated with the $i$-th power of $\mu^{-2} g$.} Algebraic properties of this element, i.e., its permutation relations with the algebra generators are described in the proposition \ref{proposition4.11}. 
\medskip

With the orthogonal $O(k)$ type Yang-Baxter matrix $R$ one associates 
yet another rank-one projector $\rho_R(a^{(k)})$. The corresponding group-like element in the
FRT-algebra ${\cal M}(R,P)$
is the $k$-th order polynomial in the algebra generators
called the {\em quantum determinant of} 
$M$
\be
\lb{qdetT}
{ \rdet_q} M := \tr_{(1,2,\dots k)} (M_{\overline{1}}\,M_{\overline{2}}\dots M_{\overline{k}}
\,\rho_R(a^{(k)}))\,=\, q^{k(k-1)}\, \goth{e}_k
\ee
(for the latter equality see (\ref{qdet-1})). 
By a slight abuse of notation we will call element $\goth{e}_k$, in the $O(k)$ type QM-algebra 
${\cal M}(R,F)$, the  quantum determinant of
$M$.
\smallskip

%
In the rest of this subsection
we derive the permutation relations between the quantum
 determinant $\goth{e}_k$  and the components of the quantum
matrix $M$. Our method of derivation does not depend on the choice of the Yang-Baxter matrix $F$, so we carry out the consideration for a general orthogonal QM-algebra.

\begin{prop}\lb{proposition4.16a}
Let ${\cal M}(R,F)$ be the QM-algebra of the orthogonal $O(k)$ type. 
Define two operators $O,O^{-1}\in {\rm End}(V)$:
\be O_1\, :=\, {k_q (1+q^{2-k})\over(q+q^{1-k})}\, \tr_{\!(2,\dots ,k+1)}
\left(\rho_R(a^{(k)\uparrow 1})\, F^{-1}_1 F^{-1}_2\dots F^{-1}_k\right)\,
\lb{O}\ee
and
\be
\lb{O-}
O^{-1}_1\,  =\, {k_q (1+q^{2-k})\over(q+q^{1-k})}\, \tr_{\!(2,\dots ,k+1)}
\left( F_k\dots F_2 F_1\, \rho_R(a^{(k)\uparrow 1})\right)\, .
\ee

\smallskip
Then

\medskip\noindent
{\rm (a)} the commutation relations of the quantum determinant with the quantum
matrix elements are given by
\be
\lb{Ma_k}
M\, \goth{e}_k\, =\, \goth{e}_k\, ( O^{-1} M O)\, ;\ee

\noindent
{\rm (b)} the operators $O$ and $O^{-1}$ are inverse to one another.
\end{prop}

\nin {\bf Proof.~} (a)~ A check of the commutation relation (\ref{Ma_k}) is similar to
the proof of formula (\ref{g-perm}) (see lemma 4.13 in \cite{OP-BMW}), although a bit 
more sophisticated. For the calculation we use  a version of the ``vanishing" condition in 
(\ref{spec4}), involving the operator
$\rho_R(a^{(k)\uparrow 1}a^{(k)}a^{(k)\uparrow 1})$. To derive it we use the next 
Lemma.

\begin{lem} Let $R$ be the Yang--Baxter matrix  of the BMW type. Assume that the 
conditions
(\ref{q-mu-restrict-2}) for the parameters $q$ and $\mu$ are satisfied. Then
\be\begin{array}{l} j_q^2\,\rho_R(a^{(j)\uparrow 1}a^{(j)}a^{(j)\uparrow 1})=
\rho_R(a^{(j)\uparrow 1})+q^j(j-1)_q\,\rho_R(a^{(j)\uparrow 1}\sigma_1^-(q^{-2j})
a^{(j)\uparrow 1})\\[1em]
\ \ \ \ \ \ \ \ \ \ \ \ \ \ \ \ \ \ \ \ \ \ \ \ \ \ \
+{\displaystyle \frac{\mu^{-1}q^{-2j+2}
(q-q^{-1})^2(j-1)_q\, j_q}{(\mu^{-1}q^{1-2j}+1)(\mu^{-1}q^{3-2j}+1)}}\,
\rho_R(a^{(j)\uparrow 1}\kappa_1a^{(j)\uparrow 1})\ \lb{aj1ajaj1}\end{array}\ee
for all $j=1,\dots ,k$.
\end{lem}

\nin {\bf Proof.~} Using the second expression for $a^{(j)}$ in eq.(\ref{a^k})
and the relations (see \cite{OP-BMW}, proposition 2.2)
\be
\lb{aj-am}
a^{(j)}a^{(j-1)}=a^{(j-1)}a^{(j)}=a^{(j)}\, ,
\ee 
we find
\be \rho_R(a^{(j)\uparrow 1}a^{(j)}a^{(j)\uparrow 1})=\frac{q^{j-1}}{j_q}\,
\rho_R\Bigl( a^{(j)\uparrow 1}(1+\frac{q^{2-2j}-1}{q-q^{-1}}\,\sigma_1+
\frac{q^{2-2j}-1}{\mu^{-1}q^{3-2j}+1}\,\kappa_1) a^{(j)\uparrow 1}\Bigr)\ .\lb{aj1p}\ee
On the other hand, by eq.(\ref{baxterization}) and the idempotency of $a^{(j)}$,
\be\begin{array}{l} \rho_R (a^{(j)\uparrow 1}\sigma_1 a^{(j)\uparrow 1})=
{\displaystyle \frac{q^j}{j_q}}\,\rho_R (a^{(j)\uparrow 1})-
{\displaystyle \frac{q^j}{j_q}}\,
\rho_R(a^{(j)\uparrow 1}\,\sigma_1^-(q^{-2j})a^{(j)\uparrow 1})\\[1em]
\ \ \ \ \ \ \ \ \ \ \ \ \ \ \ \ \ \ \ \ \ \ \ \
-{\displaystyle \frac{q-q^{-1}}{\mu^{-1}q^{1-2j}+1}}\,
\rho_R(a^{(j)\uparrow 1}\kappa_1a^{(j)\uparrow 1})\ .\end{array}\lb{aj1v}\ee
Substituting the expression for the operator $\rho_R (a^{(j)\uparrow 1}
\sigma_1 a^{(j)\uparrow 1})$ from eq.(\ref{aj1v}) into eq.(\ref{aj1p}), we arrive at eq.
(\ref{aj1ajaj1}).\hfill$\blacktriangleright$

\medskip
Assume that the Yang--Baxter matrix  $R$ is of the type $O(k)$. Then, taking into 
account the conditions $\mu=q^{1-k}$ and  (\ref{spec4}),
the identity (\ref{aj1ajaj1}) for $j=k$ specifies to
\be\lb{vanish-2}
\rho_R(a^{(k)\uparrow 1})\, =\,
k_q^2\, \rho_R(a^{(k)\uparrow 1}a^{(k)}a^{(k)\uparrow 1})\, -\,
{q^{1-k} (q-q^{-1})^2 k_q (k-1)_q\over (1+q^{-k})(1+q^{2-k})}\,
\rho_R(a^{(k)\uparrow 1}\, \kappa_1\, a^{(k)\uparrow 1})\, .\ee
This is the needed version of the vanishing condition.

\medskip
Now we calculate
\ba
\nonumber
M_1\left( \goth{e}_k\, \rho_R(a^{(k)\uparrow 1})\right)\, =\,
M_1(F_1 \dots F_k)\left( \goth{e}_k\, \rho_R(a^{(k)})\right)(F_1 \dots F_k)^{-1}\, =\, q^{k(1-k)}\, 
M_{\overline{1}} M_{\overline{2}}\dots
M_{\overline{k+1}}\ \rho_R(a^{(k)\uparrow 1})\hspace{-0.5cm}
\\[2mm]
\lb{upc1}
\, =\, k_q^2\, \goth{e}_k\,  \rho_R(a^{(k)\uparrow 1}a^{(k)})M_{\overline{k+1}}\ 
\rho_R(a^{(k)\uparrow 1})
- {\textstyle {q^{1-k} (q-q^{-1})^2 k_q (k-1)_q\over (1+q^{-k})(1+q^{2-k})}}\,
\rho_R(a^{(k)\uparrow 1}\kappa_1)\,M_1\, \rho_R(a^{(k)\uparrow 1}) \goth{e}_k .
\hspace{1.0cm}
\ea
Here in the second equality we used the relation (\ref{spec-aM}) and the 
``shifting" relation (see, e.g., \cite{OP-BMW}, lemma 4.3)
$$
F_1\, F_{2}\dots F_k\cdot M_{\overline 1}\, M_{\overline {2}}\dots
M_{\overline k}\, = \, M_{\overline 2}\, M_{\overline {3}}\dots
M_{\overline{k+1}}\cdot F_1\, F_{2}\dots F_k\, .
$$
In the third equality in the computation (\ref{upc1})
we used the condition~ (\ref{vanish-2}),~
then applied permutation relations~  (\ref{qmai})~ and 
$R_i M_{\overline j} = M_{\overline j} R_i\;\; \forall\, j\neq i,i+1$,~ and then used again 
the relation~(\ref{spec-aM}).

\smallskip
We continue the calculation by taking the traces in the spaces from 2 to $k+1$
in the  equality (\ref{upc1}). Due to the rank 1 property of the operator
$\rho_R(a^{(k)})$,  $\tr_{\!(2,\dots ,k+1)}\bigl(\rho_R(a^{(k)\uparrow 1})\bigr) = I_1$, 
and we find
\be\begin{array}{l} M_1\, \goth{e}_k\, =\, k_q^2\, \goth{e}_k\,\tr_{\!(2,\dots ,k+1)}\left(
\rho_R(a^{(k)\uparrow 1}a^{(k)})M_{\overline{k+1}}\right)\\[1em]
\ \ \ \ \ \ \ \ \ \ \ \ \
-\,{\displaystyle\frac{q^{1-k}(q-q^{-1})^2k_q(k-1)_q}{(1+q^{-k})(1+q^{2-k})}}\,
\tr_{\!(2,\dots ,k+1)}\left( K_1\,\rho_R(a^{(k)\uparrow 1})\right)\, M_1\, \goth{e}_k\, .
\end{array}\lb{mak0}\ee
The traces in the last term in the left hand side 
of eq.(\ref{mak0}) can be calculated with the help of the
identity
\be
\tr_{\!(2,\dots ,k)}\left(\rho_R(a^{(k)})\right)\, =\,
q^{k(k-1)}\, {D_R}_1\,\Tr{2,\dots ,k}\left(\rho_R(a^{(k)})\right)\, 
=\, {\displaystyle \frac{q^k(1+q^{-k})}{(1+q^{2-k})k_q}}\, {D_R}_1\ .
\lb{heid}
\ee
Here in the first equality we used the relation (\ref{spec-a1}); in the second  equality 
we used eq.(\ref{spec1}) ~$k-1$ times and substituted $\mu=q^{1-k}$.

\smallskip
Therefore,
\be
\tr_{\!(2,\dots ,k+1)}\left( K_1\,\rho_R(a^{(k)\uparrow 1})\right)\, =\,
\frac{q^k(1+q^{-k})}{(1+q^{2-k})k_q}\, \Tr{2}K_1=
{q(1+q^{-k})\over (1+q^{2-k}) k_q}\, I_1\,
\lb{mak2}\ee
due to the relation (\ref{traceDK}) with $\mu =q^{1-k}$.

\smallskip
After substituting the expression (\ref{mak2}) and rearranging 
similar terms, eq.(\ref{mak0}) simplifies to
\be {(q+q^{1-k})^2\over (1+q^{2-k})^2}\, M_1\, \goth{e}_k\, =\,
k_q^2\, _k\,\tr_{\!(2,\dots ,k+1)}\left(
\rho_R(a^{(k)\uparrow 1}a^{(k)})M_{\overline{k+1}}\right) .\lb{mak01}\ee

It is left to evaluate the traces in the right hand side of eq.(\ref{mak01}).

\smallskip
Due to the rank 1 property of the operator $\rho_R(a^{(k)})$, we have
\be \rho_R(a^{(k)\uparrow 1})\, X_{1,\dots ,k+1}\,\rho_R(a^{(k)\uparrow 1})=
\tr_{\!(2,\dots ,k+1)}\Bigl( \rho_R(a^{(k)\uparrow 1})\, X_{1,\dots ,k+1}\Bigr)\
\rho_R(a^{(k)\uparrow 1})\ \lb{x1kp1a}\ee
for any operator $X\in {\rm End}(V^{\otimes (k+1)})$.

\smallskip
Let $Y$ be an operator acting in the space $V$. Then
\be\begin{array}{l}
\tr_{\!(2,\dots ,k+1)}\Bigl(
\rho_R(a^{(k)\uparrow 1}\, a^{(k)})\, Y_{\overline{k+1}}\Bigr)\,
=\, \tr_{\!(2,\dots ,k+1)}\Bigl(
\rho_R(a^{(k)})\, Y_{\overline{k+1}}\,\rho_R(a^{(k)\uparrow 1})\Bigr)\, \\[1em]
\ \ \ \ \ =\, \tr_{\!(2,\dots ,k+1)}\Bigl(
\rho_R(a^{(k)})\cdot F_k\dots F_2\, F_1\cdot Y_1\cdot F_1^{-1}\, 
F_2^{-1}\dots F_k^{-1}\cdot
\rho_R(a^{(k)\uparrow 1})\Bigr)\, \\[1em]
\ \ \ \ \ =\, \tr_{\!(2,\dots ,k+1)}\Bigl(
{}F_k\dots F_2\, F_1\cdot Y_1\, \rho_R(a^{(k)\uparrow 1})\cdot
{}F_1^{-1}\, F_2^{-1}\dots F_k^{-1}\cdot\rho_R(a^{(k)\uparrow 1})\Bigr)\, \\[1em]
\ \ \ \ \ =\, {\displaystyle \frac{q+q^{1-k}}{k_q(1+q^{2-k})}}\
\tr_{\!(2,\dots ,k+1)}\Bigl(
{}F_k\dots F_2\, F_1\cdot Y_1\, O_1\,\rho_R(a^{(k)\uparrow 1})\Bigr)\,\\[1em]
\ \ \ \ \ =\, {\displaystyle \frac{q+q^{1-k}}{k_q(1+q^{2-k})}}\
\tr_{\!(2,\dots ,k+1)}\Bigl(
\rho_R(a^{(k)\uparrow 1})\cdot F_k\dots F_2\, F_1\cdot 
\rho_R(a^{(k)\uparrow 1}) \Bigr)\
Y_1\, O_1\, \\[1em]
\ \ \ \ \ =\, {\displaystyle \frac{(q+q^{1-k})^2}{k_q^2(1+q^{2-k})^2}}\
\tr_{\!(2,\dots ,k+1)}\Bigl( \rho_R(a^{(k)\uparrow 1}) \Bigr)\
O_1^{-1}\, Y_1\, O_1\, =\,
{\displaystyle \frac{(q+q^{1-k})^2}{k_q^2(1+q^{2-k})^2}}\
O_1^{-1}\, Y_1\, O_1\ .
\end{array}\lb{omyo}
\ee
Here in the first equality we moved the operator $\rho_R(a^{(k)\uparrow 1})$
cyclically to the right under the trace; in the second equality we substituted 
the definition of $Y_{\overline{k+1}}$; in the third
equality we used the twist relations (\ref{sovm}) to move the operator
$\rho_R(a^{(k)})$ rightwards; in the fourth equality we used eq.(\ref{x1kp1a})
with $X_{1,\dots ,k+1}=F_1^{-1}\, F_2^{-1}\dots F_k^{-1}$ and the definition
(\ref{O}) of the operator $O$; in the fifth equality we replaced $\rho_R(a^{(k)\uparrow 1})$
by $\rho_R(a^{(k)\uparrow 1})\rho_R(a^{(k)\uparrow 1})$ 
using the idempotency
of the operator $\rho_R(a^{(k)\uparrow 1})$,
moved one of the
operators $\rho_R(a^{(k)\uparrow 1})$ cyclically to the left under the trace
and moved $Y_1O_1$ out of the trace; in the sixth equality we used eq.(\ref{x1kp1a})
with $X_{1,\dots ,k+1}=F_k\dots F_2F_1$ and the definition (\ref{O-}) of the
operator $O^{-1}$; in the seventh equality the rank 1 property
of the operator $\rho_R(a^{(k)})$ was used again.

\medskip
Substituting the result of the calculation (\ref{omyo}) with $Y=M$ into
eq.(\ref{mak01}), we finally establish the relation (\ref{Ma_k}).

\medskip
(b)~ To check the relation $O^{-1} O = I$ we use the same calculation (\ref{omyo}),
but now with $Y=I$: in this case
\be
\begin{array}{l}
\tr_{\!(2,\dots ,k+1)}\Bigl(
\rho_R(a^{(k)\uparrow 1}\, a^{(k)})\, Y_{\overline{k+1}}\Bigr)\, \rule[-9pt]{0.5pt}{20pt} 
_{_{\, Y=I}}
\, =\,
\tr_{\!(2,\dots ,k+1)}\Bigl(
\rho_R(a^{(k)\uparrow 1}\, a^{(k)})\Bigr)\, \\[1em]
\ \ \ \ \ =\,
{\displaystyle\frac{q^{k-1}}{k_q}}\
\tr_{\!(2,\dots ,k+1)}\Bigl(
\rho_R(\sigma_1(q^{2-2k})\, a^{(k)\uparrow 1})\, \Bigr)\,
 =\,
{\displaystyle\frac{q^{-1}(1+q^{-k})}{k_q^2(1+q^{2-k})}}\
\Tr{2}\Bigl(\rho_R(\sigma_1(q^{2-2k}))\Bigr)\,\\[1em]
\ \ \ \ \ =\,
{\displaystyle \frac{(q+q^{1-k})^2}{k_q^2(1+q^{2-k})^2}}\ I_1\ .
\end{array}
\lb{omyo2}
\ee
Here in the second equality we used the second expression for $a^{(k)}$ in the 
definition
(\ref{a^k}), the cyclic property of the trace and the identities (\ref{aj-am})
for the $q$-antisymmetrizers; in the third equality we used the 
calculation (\ref{heid}); in
the fourth equality we evaluated the R-traces with the help of the relations
(\ref{trC-D}), (\ref{traceDK}) and (\ref{traceD}).

We are done by comparing the results of calculations (\ref{omyo2}) and (\ref{omyo}) for $Y=I$. 
\hfill$\blacksquare$

\begin{rem}{\rm
Using a similar techniques, one can verify independently
relation $OO^{-1} = I$ which is of course equivalent to
 $O^{-1}O = I$ for a finite-dimensional $V$. This is done by a calculation, 
 parallel to (\ref{omyo}),
where one uses the copies $Y_{\underline{j}}$, defined in eq.(\ref{kopii-2})$\,$,
instead of the copies $Y_{\overline{j}}$. \hfill$\blacktriangleright$
}\end{rem}

\begin{rem}{\rm For the compatible pair $\{ R,R\}$ of Yang--Baxter matrices, we can
evaluate all traces in the expression (\ref{O-}) moving cyclically the operator
$\rho_R(a^{(k)\uparrow 1})$ under the trace and using characteristic 
properties of the $q$-antisymmetrizers $a^{(n)}\sigma_i =
\sigma_i a^{(n)}=-q^{-1}a^{(n)}$ $\forall 1\leq i<n$ (see \cite{OP-BMW}, 
proposition 2.2)
and (\ref{heid}). Thus we obtain $O^{-1}\!=O=(-1)^{k-1}\, I$. Likewise, the 
operator $G$ (\ref{G})
trivializes for the compatible pair $\{ R,R\}$: we have $G=I$. 
		
Scalar forms of the operators $O$ and $G$ in the reflection equation 
algebra ${\cal M}(R,R)$ was expected since the elements $\goth{e}_k$ and $g$ both 
belong to the center of the algebra
 (see remark 	\ref{rem3}). \hfill$\blacktriangleright$
}\end{rem}

\subsection{Representations of quantum matrix algebras}\lb{repqma}

In the next subsection we will use two representations of the  QM-algebra 
${\cal M}(R,F)$
on the same underlying vector space $V$. In the present subsection these     
representations -- we denote them by $\alpha_+$ and $\alpha_-$ --
are introduced and 
investigated.
We use the label 0 for the copy of the
vector space $V$, on which $\alpha_\pm$ are defined. Here the symbol 0 carries 
an auxiliary role; 0 does not mean a number preceding 1; there is no matrix $M_0$. The space 
$V_0$ is just another space in which our operators act.

\begin{prop}\lb{lemmsol} Let ${\cal M}(R,F)$ be the QM-algebra of the general type (no 
additional  conditions on the compatible pair $\{ R,F\}$ of Yang--Baxter matrices).
	
\smallskip
The following formulas:
\be 
\alpha_\varepsilon(M_1)_0=R_{01}^{\varepsilon}\, F_{01}\ ,\ {\mbox or,\, shortly,}\ R_0^\varepsilon F_0\lb{solm}
\ee
where $\varepsilon$ is the sign, $\varepsilon=\pm 1$, define two representations $\alpha_\varepsilon$ of the
	 QM-algebra ${\cal M}(R,F)$ on the vector space $V$.
\end{prop}

\noindent {\bf Proof.}~ A straightforward calculation, with the help of the Yang-Baxter
equations (\ref{YB}) and the twist relations (\ref{sovm}) for the compatible pair 
$\{ R,F\}$ of Yang--Baxter matrices.\footnote{The solution $\alpha_+$  
was noted in \cite{IOP2}. The solution $\alpha_-$ becomes 
trivial in the case
of the reflection equation algebra, when $F=R$.} \hfill$\blacksquare$

\begin{rem}{\rm There are many ways to construct representations of the  QM-algebra 
${\cal M}(R,F)$
using the compatible Yang--Baxter matrix  pair $\{R,F\}$.		
In particular, for ${\cal M}(R,F)$ of the general type, there is
yet another representation $\beta$ on the vector space $V$. It is given by
\[ \hspace{6.5cm}\beta(M_1)_0= F_{01}^2\ ,\ {\mbox or}\ F_0^2
\mbox{\hspace{6.6cm}$\blacktriangleright$}\]
}\end{rem}

\medskip
The quantum matrix string $M_{\overline{1}}\, M_{\overline{2}}\dots M_{\overline{j}}\, $, 
which we always use to construct elements of the characteristic subalgebra 
${\cal C}(R,F)$ 
acquires a quite simple expression in the representations $\alpha_\pm$.

\begin{lem}\lb{lemmastrFR}
We have
\be \alpha_\varepsilon(M_{\overline{1}}M_{\overline{2}}\dots M_{\overline{j}})_0=
R^\varepsilon_0R^\varepsilon_1\dots R^\varepsilon_{j-1}
 F_{j-1}\dots F_1F_0
 \lb{strFRrep}
 \ee
for any $j>0$.
\end{lem}

\noindent {\bf Proof.}~ Induction on $j$:
$$
\begin{array}{l} 
\alpha_\varepsilon(M_{\overline{1}}M_{\overline{2}}\dots M_{\overline{j+1}})_0=
(R^\varepsilon_0R^\varepsilon_1\dots R^\varepsilon_{j-1} F_{j-1}\dots 
F_1F_0) F_{j}\dots F_2F_1\,
R^\varepsilon_0F_0\, (F_{j}\dots F_2F_1)^{-1}
\\[5pt]
\ \ \ \ \ \ \ \ \ \
=R^\varepsilon_0R^\varepsilon_1\dots R^\varepsilon_{j-1} 
(F_{j-1}F_j)(F_{j-2}F_{j-1})\dots (F_0F_1)\,
R^\varepsilon_0F_0\, (F_{j}\dots F_2F_1)^{-1}
\\[5pt]
\ \ \ \ \ \ \ \ \ \
=R^\varepsilon_0R^\varepsilon_1\dots R^{\varepsilon}_{j-1}\, R^{\varepsilon}_jF_j\, (F_{j-1}F_j)(F_{j-2}F_{j-1})\dots (F_0F_1) (F_{j}\dots F_2F_1)^{-1}
\\[5pt]
\ \ \ \ \ \ \ \ \ \
=R^\varepsilon_0R^\varepsilon_1\dots R^\varepsilon_{j} F_j\dots 
F_1F_0\ .
\end{array}
$$
In the first equality we used the induction assumption and the definition
of the matrix $M_{\overline{j+1}}$. In the second equality we shifted the
elements of the string $F_{j-1}\dots F_1F_0$ as far as possible to the right into the string
$F_{j}\dots F_2F_1$; as a result they formed pairs. 
In the third equality we used the 
compatibility relations for the $\{R,F\}$ pair to move the product 
$R^\varepsilon_0F_0$
leftwards. In the fourth equality we reassemled the product of pairs of the
operators $F$ back into the product of two
strings $F_{j-1}\dots F_1F_0\cdot F_{j}\dots F_2F_1$. The proof is 
finished.\hfill$\blacksquare$

\medskip
We shall now calculate the image of the contraction $g$ in the representations
$\alpha_\pm$.

\begin{prop}\lb{lemimageg}
For QM-algebras of the BMW type the images of the contraction $g$ in the 
representations
$\alpha_\pm$ are proportional to the operator $G^{-1}$,
\be \alpha_\pm(g)=\mu^2 G^{-1}\ .\lb{imageg}\ee
\end{prop}

\noindent {\bf Proof.}~ The proof is a direct calculation:
$$
\begin{array}{l} \alpha_\pm(g)_0=\eta^{-1}\,\Tr{1,2}
(K_1R^\varepsilon_0R^\varepsilon_1F_1F_0)=
\eta^{-1}\mu^2\,\tr_{\!(1,2)}(K_1K_0F_1F_0)
\\[5pt]
\ \ \ \ \ \ \ \ \ \ =\eta^{-1}\mu^2\,\tr_{\!(1,2)}(K_0F_1F_0K_1)=\mu^2\,\tr_{\!(1,2)}
(K_0F_1F_0)=\mu^2G_0^{-1}.
\end{array}
$$
Here in the first equality we used the definition (\ref{tau}) of the contraction
$g$ and expressed the string $M_{\overline{1}}M_{\overline{2}}\, $ as in the lemma
\ref{lemmastrFR}; in the second equality we used relations (\ref{bmw3}) and 
(\ref{KDD}); in the
third equality we moved $K_1$ cyclically to the right under the trace; in the
fourth equality we used the twist relation and the identity $K^2=\eta K$;
the fifth equality reproduces the formula (\ref{G}).
\hfill$\blacksquare$

\begin{prop}\lb{propimagedet}
For orthogonal QM-algebras of the type $O(k)$ the images of the quantum 
determinant 
$\goth{e}_k\in {\cal C}(R,F)$ in the representations $\alpha_+$ and $\alpha_-$ are given 
by
\be \alpha_\pm(\goth{e}_k)\,=\, (-1)^{k-1}q^{k(1-k)}\ O^{-1}\ .\lb{aplde}\ee
\end{prop}

\noindent {\bf Proof.}~ Evaluating formula (\ref{spec-aM}) in the representation 
$\alpha_+$ we find, with the help
of the relation (\ref{strFRrep}), that
$$
\rho_R(a^{(k)})\, R_{0} R_{1}\dots R_{k-1} F_{k-1}\dots F_1 F_0\, =\, q^{k(k-1)}\ 
\rho_R(a^{(k)})\ (\alpha_+(\goth{e}_k))_0\ .
$$
Multiplying both sides of this equality by the operator $F_0^{-1}$ from the right, by
the operator ${C_R}_0$ (\ref{C-D}) from the left, taking the trace in the 
 the space $V_0$ and 
using eq.(\ref{trC-D}), we find
$$  
\rho_R(a^{(k)})\, F_{k-1}\dots F_2 F_1\, =\, (-1)^{k-1}\, q^{(k+1)(k-1)}
\rho_R(a^{(k)})\ \tr_{(0)}({C_R}_0\, (\alpha_+(\goth{e}_k))_0\, F_0^{-1})\ ,
$$
because $\rho_R(a^{(k)})\, R_{j}=(-q)^{-1}\rho_R(a^{(k)})$ for all $j=1,\dots ,k-1$.

\smallskip
Taking the traces in the spaces from 2 to $k$ and using eq.(\ref{heid}), we obtain, 
upon the replacement $i\mapsto i+1$ of the labels of all the spaces, except $V_0$:
\be\begin{array}{l}
\tr_{\!(3,\dots ,k+1)}\left( \rho_R(a^{(k)\uparrow 1}) F_k\dots F_3 F_2
\right)\,\\[5pt]
\ \ \ \ \ \ \ \ \ \ \ \ \ \ \ \ \ \ \ \ =\, (-1)^{k-1}\,
{\displaystyle\frac{q^{k^2+k-1}(1+q^{-k})}{k_q(1+q^{2-k})}}\ {D_R}_2
\tr_{(0)}({C_R}_0\, (\alpha_+(\goth{e}_k))_0\, F_{02}^{-1}) .
\end{array}
\lb{eps3}
\ee

By the lemma 3.3, \cite{OP-BMW} $(\Psi_{F^{-1}})_{21}=(D_R^{-1})_1F_{12}{D_R}_2$.
Multiplying the relation (\ref{eps3}) from the right by a combination 
$(D_R^{-1})_1F_{12}$, we
take the trace in the space 2, moving the operator ${D_R}_2$ cyclically to the
right under the trace in the right hand side and applying eq.(\ref{s-inv}) for the 
Yang--Baxter matrix  $F^{-1}$. All traces in the right hand
side can now be evaluated and we find
\be 
(D^{-1}_R)_1\, \tr_{\!(2,\dots ,k+1)}\left( \rho_R(a^{(k)\uparrow 1})\,
{}F_k\dots F_2 F_1\right)\, =\, (-1)^{k-1}\,
{\displaystyle\frac{q^{k^2+k-1}(1+q^{-k})}{k_q(1+q^{2-k})}}\
{C_R}_1\, (\alpha_+(\goth{e}_k))_1\ .
\lb{eps4}
\ee
Using the definition (\ref{O-}) of the matrix $O^{-1}$ and the relation (\ref{Cinv-D}), we
establish the relation (\ref{aplde}) for the representation $\alpha_+$.

\smallskip
The calculations for the representation $\alpha_-$ are similar, we do not
repeat them.\hfill$\blacksquare$
\medskip

In the end of the subsection we apply the representation $\alpha_+$ to clarify 
properties of the operator $O$.
\begin{prop}
\lb{proposition4.16b} 
Let $\{R,F\}$ be a compatible Yang--Baxter matrix  pair, where $R$ is $O(k)$ 
type. The corresponding matrices $O$, $G$, $D_R$ satisfy following relations
\ba O_1\, F^\varepsilon_{12}\, & =& \, F^\varepsilon_{12}\, O_2\, ,\qquad 
\mbox{where $\varepsilon=\pm 1$,}
\lb{ofo-1}
\\[2pt]
R_{12}\, O_1\, O_2& =& O_1\, O_2\, R_{12}\, ,
\lb{oor-1}
\\[2pt]
[ D_R,O ]
& = &  0\,\, =\,\, [ G,O ]\, ,
\lb{commdro-1}
\\[2pt]
O^2& =& G^k\, .
\lb{o2gk-1}
\ea
\end{prop}

\noindent {\bf Proof.}~
{}For an arbitrary integer $j>0$ and an operator $X=X_{1,\dots ,j}
\in {\rm End}(V^{\otimes j})$, let
$\tilde{X}\in {\rm End}(V)$ be an operator defined by
$\tilde{X}_1=\Tr{2,\dots ,j}(X_{1,\dots ,j})$. It follows then from lemma 3.2,  
\cite{OP-BMW}  that
$$
\Tr{3,\dots ,j+1}\Bigl( F_2^{-\varepsilon} F_3^{-\varepsilon}\dots F_j^{-\varepsilon}\,
X_{1,\dots ,j}\, F_j^\varepsilon\dots F_3^\varepsilon F_2^\varepsilon\Bigr)\, =\,
\tilde{X}_1 I_2\ 
$$
for $\varepsilon=\pm 1$.

\smallskip
Therefore, for the operator
$X_{1,\dots ,k+1}=\rho_R(a^{(k)\uparrow 1})\cdot F_k\dots F_2\, F_1$, we have
\be\begin{array}{l}
\tilde{X}_1 F^\varepsilon_1\, =\, \Tr{3,\dots ,k+2}\Bigl( F_2^{-\varepsilon} 
F_3^{-\varepsilon}\dots
{}F_{k+1}^{-\varepsilon}\,\rho_R(a^{(k)\uparrow 1})\, F_k\dots F_2 F_1\,
{}F_{k+1}^\varepsilon\dots F_2^\varepsilon F_1^\varepsilon\Bigr)
\\[5pt]
\ \ \ \ \ =\, \Tr{3,\dots ,k+2}\Bigl( F_2^{-\varepsilon} F_3^{-\varepsilon}\dots
{}F_{k+1}^{-\varepsilon}\,\rho_R(a^{(k)\uparrow 1})\, F_{k+1}^\varepsilon
\dots F_2^\varepsilon F_1^\varepsilon
\, F_{k+1}\dots F_3 F_2\Bigr)
\\[5pt]
\ \ \ \ \ =\,  F_1^\varepsilon\ \Tr{3,\dots ,k+2}\Bigl(\,\rho_R(a^{(k)\uparrow 2})
\, F_{k+1}\dots F_3 F_2\Bigr)\, =\, F_1^\varepsilon \tilde{X}_2\ .
\end{array}\label{tildeX-F}\ee
Here we used the 
compatibility relations for the $\{R,F\}$ pair.

\smallskip
Noticing that for the chosen $X$ the operator $\tilde{X}$ is proportional to the operator 
$O^{-1}$,
(see the definition (\ref{O-}) and the relation (\ref{spec-a1})$\,$) we see that the 
calculation  (\ref{tildeX-F}) proves the relation (\ref{ofo-1}).

\medskip
The commutation relation (\ref{Ma_k}) between the quantum matrix and the quantum 
determinant, in
the representation $\alpha_+$, reads, by the proposition \ref{propimagedet},
\be R_0\, F_0\, O_0^{-1}\, =\, O_0^{-1}\, O_1^{-1}\, R_0\, F_0\, O_1\ .
\lb{oor1-1}
\ee
Applying the relation (\ref{ofo-1}) to eq.(\ref{oor1-1}), we establish the relation
(\ref{oor-1}).

\medskip
By eq.(\ref{FCC}) the following identity
\be\begin{array}{l} F_j\dots F_2 F_1\, {D_R}_1\, =\, F_j\dots F_2\, F_1\,
({D_R}_1\, {D_R}_2\dots {D_R}_{j+1})\, ({D_R^{-1}}_2\, {D_R^{-1}}_3\dots 
{D_R^{-1}}_{j+1})
\\[1em]
\ \ \ \ \ \ \ \ \ \ =\, ({D_R}_1 {D_R}_2\dots {D_R}_{j+1})
\, F_j\dots F_2 F_1\, ({D_R^{-1}}_2 {D_R^{-1}}_3\dots {D_R^{-1}}_{j+1})
\ \end{array}\lb{od1-1}\ee
holds for any $j>0$.

\smallskip
The cyclic properties of the trace, together with eq.(\ref{od1-1}), shows that
the expression (\ref{O-}) commutes with ${D_R}_1$, which proves the left equality in
 (\ref{commdro-1}).

\medskip
The right equality in (\ref{commdro-1}) is 
a consequence of the commutativity of the characteristic subalgebra ${\cal C}(R,F)$ 
because
the images of the elements $g,\goth{e}_k\in {\cal C}(R,F)$ in the representations 
$\alpha_\pm$
are given by the formulas (\ref{imageg}) and (\ref{aplde}).

\medskip
The relation (\ref{o2gk-1}) is the matrix image of the reciprocal relation 
(\ref{reciprocal})
with $i=0$, $g^k=\goth{e}_k^2$, in the representation $\alpha_+$.

\smallskip
The proof is finished.\hfill$\blacksquare$

\subsection{Resolution of reciprocal relations}\lb{components}

In this subsection we consider the extended QM-algebra ${\cal M^{^\bullet\!}}(R,F)$ of 
the orthogonal type $O(k)$. 
We shall use shorthand notation $O_k(R,F)$ for these algebras from now on.

We employ the invertibility of the
contraction $g$ to resolve the reciprocal relations with respect to the first $\lfloor k/2\rfloor$ 
elementary sums.

The case $k=2\ell$ requires a more detailed investigation. For $k=2\ell$ we 
can resolve the reciprocal relations under the 
assumption that the element $g^{-\ell}\goth{e}_{2\ell}$ is central.
To motivate this assumption we first notice that $(g^{-\ell}\goth{e}_{2\ell})^2=1$
by eq.(\ref{reciprocal}) in which $k=2\ell$ and $i=0$.
Next, the following proposition holds.
\begin{prop}\lb{proposition4.19}
Assume, for the (extended) QM-algebra $O_{2\ell}(R,F)$, that at least  one of the 
matrices $\alpha_\pm(\goth{e}_\ell )$ is non-degenerate.
Then the element $g^{-\ell}\goth{e}_{2\ell} $ is central in $O_{2\ell}(R,F)$.
\end{prop}
\begin{rem}{\rm
In the classical
case ($R=F=P$), both matrices $\alpha_\pm(\goth{e}_\ell )$ are proportional
to the identity matrix, the coefficient of proportionality equals
${(2\ell )!}/{(\ell !)^2}$~ (the rank of the $\ell$-th
$q$-antisymmetrizer), so the assumption of Proposition \ref{proposition4.19}
holds on the perturbation level.}
\hfill$\blacktriangleright$
\end{rem}

\noindent {\bf Proof.}~
By the propositions \ref{proposition4.11} and \ref{proposition4.16a} 
the element $g^{-\ell} \goth{e}_{2\ell}$ commutes with the components of the quantum matrix 
$M$ 
if and only if  the matrix $G^{\ell} O^{-1}$ is scalar.

Taking the images of the reciprocal relations (\ref{reciprocal}) in the representations
$\alpha_{\varepsilon}$, $\varepsilon=\pm 1$ we find, in particular, that 
$\alpha_\varepsilon(g)^\ell \alpha_\varepsilon(\goth{e}_\ell )=
\alpha_\varepsilon(\goth{e}_{2\ell}) \alpha_\varepsilon(\goth{e}_\ell)$ which, by the propositions 
\ref{lemimageg}, \ref{propimagedet} and under the assumption of the present
proposition, implies that $G^{-\ell}=-O^{-1}$. 
It follows then  that the
element $g^{-\ell} \goth{e}_{2\ell}$ is central. 
\hfill$\blacksquare$

\begin{defin}
	\lb{comp-2ell}
Assume that the element $g^{-\ell}\goth{e}_{2\ell}$ belongs to the center of the  QM-algebra 
$O_{2\ell}(R,F)$. Then, its quotient algebra by the relation
\be
\lb{quotients-O}
g^{-\ell}\goth{e}_{2\ell} \, =\, 1, \ \ \mbox{respectively}\ \  g^{-\ell}\goth{e}_{2\ell} \, =\,  -1,
\ee
is called the component $O^{+}_{2\ell}(R,F)$, respectively $O^{-}_{2\ell}(R,F)$.
\end{defin}

In fact we defined the decomposition of the algebra $O_{2\ell}(R,F)$ into the direct sum 
of its two {\em components} $O^{\pm}_{2\ell}(R,F)$. Note 
that classically, that is, for $R=F=P$ and under the restriction $g=1$, $O^{\pm}_{2\ell}$ 
become the algebras of functions on two connected components of the orthogonal 
group $O_{2\ell}$. 

Now we are ready to resolve the reciprocal relations. 

\begin{theor}
\lb{theorem4.23}
In the assumption of the definition \ref{comp-2ell}, the 
reciprocal relations (\ref{reciprocal}) in the component algebras
$O^{\pm}_{2\ell}(R,F)$ can be resolved with respect to $\goth{e}_j$, $j\geq \ell$:
\be
\lb{reciprocal2} 
\goth{e}_{\ell +i} = \pm g^i\, \goth{e}_{\ell -i}\,\qquad \forall\; i=0,1,\dots ,\ell\, ,
\ee
respectively. Note that in the $O^-_{2\ell}$ case, this implies
that $\goth{e}_{\ell}=0$.\footnote{A classical counterpart is: if $X$ is an
orthogonal operator acting in a space $V$ of dimension $2\ell$ and 
${\mathrm{det}}\,X=-1$
then ${\mathrm{tr}}_{\wedge^\ell V}(X)=0$, where $\wedge^j V$, $j=2,3,\dots \,$, are 
the wedge tensor powers of the space $V$.}
\smallskip

{}For the $O_{2\ell -1}(R,F)$ 
the reciprocal relations
(\ref{reciprocal}) admit the resolution
\be
\lb{reciprocal3} 
\goth{e}_{\ell+i}\, =\, g^{i+1/2} \goth{e}_{\ell-1-i}\,\qquad
\forall\; i=0,1,\dots \ell-1\, ,
\ee
where we introduce 
\be
\lb{root-g}
g^{1/2} :=  g^{1-\ell} \goth{e}_{2\ell -1},
\ee
so that $g=(g^{1/2})^2$ and $\goth{e}_{2\ell -1}=(g^{1/2})^{2\ell-1}$.

\end{theor}

\noindent {\bf Proof.}~
Reciprocal relations (\ref{reciprocal}) are not independent. Since the element $g$ is invertible,
the relations are invariant under the involution $i\leftrightarrow k-i$. 
The only 
unpaired relation corresponds to $i=0$.
Resolution of this relation depends on the parity of $k$.

If $k=2\ell$, eq.(\ref{reciprocal}) for $i=0$ splits in two factors
$$
g^{2\ell}-\goth{e}_{2\ell}^2 = (g^{\ell}+\goth{e}_{2\ell})(g^{\ell}-\goth{e}_{2\ell})=0\ ,
$$ 
which implies the formulas (\ref{reciprocal2}) 
for the component algebras 
$O^{\pm}_{2\ell}(R,F)$. 

If $k=2\ell-1$, eq.(\ref{reciprocal}) for $i=0$ assumes a form
$$
g= (g^{1-\ell} \goth{e}_{2\ell-1})^2
$$
which allows one to define square root of $g$ by the formula (\ref{root-g}). 
\medskip

The check of the rest of the proposition is straightforward. \hfill$\blacksquare$

\begin{cor}\lb{cor4.24}
 Characteristic subalgebras of the extended orthogonal QM-algebras  
 $O_{2\ell-1}(R,F)$ and (in the assumption of definition \ref{comp-2ell}) of the component 
 algebras
$O^+_{2\ell}(R,F)$ and $O^-_{2\ell}(R,F)$  are generated, respectively, by the sets of 
elements $g^{1/2}\cup\{ \goth{e}_i\}_{i=0}^{\ell-1}$,~~ $g\cup\{ \goth{e}_i\}_{i=0}^\ell$ and  
$g\cup\{ \goth{e}_i\}_{i=0}^{\ell-1}$.
\end{cor}

\noindent {\bf Proof.}~
Finite height condition (\ref{spec4}) implies that for the orthogonal $O(k)$ type 
QM-algebras
$\goth{e}_{k+i}=0$ $\forall\, i\geq 1$. The corollary then follows directly from  proposition 
\ref{gen-char} and  theorem \ref{theorem4.23}.
\hfill$\blacksquare$

\section*{Acknowledgments}

The authors Pavel Saponov for fruitful discussions and valuable remarks.
The work of the second author (P. P.) was  supported
by the Basic Research Program of HSE University.
\bigskip

\appendix

\section{ Braided coproduct on  quantum matrix algebras }\lb{appendix-A}

\def\theequation{\thesection.\arabic{equation}}
\makeatletter\@addtoreset{equation}{subsection}\@addtoreset{equation}{section}\makeatother

The ${\cal M}(R,P)$ FRT-algebras possess the coproduct structure.
In this appendix we discuss its generalization, usually 
called the {\em braided coproduct}, to QM-algebras ${\cal M}(R,F)$. The braided coproduct 
structure is 
is known, mainly in the case of reflection equation algebras ${\cal M}(R,R)$, in the quantum group 
theory community. 
We describe it here for general QM-algebras.
	
Before we begin, we note that the construction of the group-like elements described in  
subsection \ref{qdetorth} for orthogonal/BMW type FRT-algebras carries over without any 
changes to general orthogonal/BMW type QM-algebras.
\medskip

(i)	 
To define the braided coproduct on ${\cal M}(R,F)$ let us interpret the coproduct on the algebra ${\cal M}(R,P)$ 
as follows. The matrices $M\otimes I$
and $I\otimes M$ element-wise commute in the vector space 
${\cal M}(R,P)\otimes{\cal M}(R,P)$ endowed with  the canonical algebra structure; the canonical
algebra structure on the tensor product of the vector spaces underlying two associative algebras 
$\cal{A}$ 
and $\cal{B}$ means exactly that any elements $x\otimes 1$ and $1\otimes y$,
$x\in\cal{A}$ and $y\in\cal{B}$, commute (while, as algebras, 
$\cal{A}\;\otimes$$\, 1$ is isomorphic to $\cal{A}$ and $1\otimes \cal{B}$ is 
isomorphic to $\cal{B}$). So the algebra 
${\cal M}(R,P)\otimes{\cal M}(R,P)$
is the same as the unital algebra
generated by two commuting quantum matrices that is 
by the elements  
$\|\dM_a^b\|_{a,b=1}^{\mbox{\footnotesize\sc n}}$ and 
$\|\ddM_a^b\|_{a,b=1}^{\mbox{\footnotesize\sc n}}$ with the defining relations	
\be\lb{twocooqmatr}R_1 
\dM_{1}\dM_{2} = \dM_{1}\dM_{2}R_1\ , \ 
R_1 \ddM_{1}\ddM_{2} = \ddM_{1}\ddM_{2}R_1\ ,\ee
and 
\be\lb{comcop}\ddM_1\dM_2=\dM_2\ddM_1\ .\ee 
The coproduct is the homomorphism
$\Delta:{\cal M}(R,P)\to {\cal M}(R,P)^{\otimes 2}$
defined on generators by 
\be\label{sopP} \Delta (M_b^a) = \sum_c 
\dM_c^a \ddM_b^c\ .\ee
We shall henceforth refer to the algebra ${\cal M}(R,P)^{\otimes 2}$ as ${\cal M}^{(2;P)}(R)$, 
to avoid the tensor sign in the name of the algebra
with a view to endow 
the space 
${\cal M}(R,F)^{\otimes 2}$ with a non-canonical algebra structure.
The upper `2' in the symbol ${\cal M}^{(2;P)}(R)$ means 
that the algebra is generated by two copies of 
the quantum matrix, while the upper `$P$' means that we use the flip operator $P$ 
to create the second 
copy 
$\dM_2=P\dM_1P$, of the quantum matrix $\dM_1$, which, according to (\ref{comcop}), 
commutes with the quantum matrix $\ddM_1$.
 
			
\vskip .4cm
(ii) Now for a general compatible pair $\{R,F\}$ define the unital algebra 
generated by the two quantum matrices ---
by the elements  $\|\dM_a^b\|_{a,b=1}^{\mbox{\footnotesize\sc n}}$ and 
$\|\ddM_a^b\|_{a,b=1}^{\mbox{\footnotesize\sc n}}$ which satisfy the relations
\be\lb{twocooqmatrq}R_1 \dM_{\overline{1}}\dM_{\overline{2}} = 
\dM_{\overline{1}}\dM_{\overline{2}}R_1\ , \ 
R_1 \ddM_{\overline{1}}\ddM_{\overline{2}} = 
\ddM_{\overline{1}}\ddM_{\overline{2}}R_1 ,\ee
but the copies $\dM$ and $\ddM$ no longer commute; instead, they have the 
following commutation relations
\be\lb{comcopF1}\ddM_{\overline{1}}\dM_{\overline{2}}=
\dM_{\overline{2}}\ddM_{\overline{1}}\ .\ee 
For a strict skew-invertible $F$ these commutation relations allow to express 
any product $\dM_b^a \ddM_d^c$ as a linear combination of 
products $\ddM_{d'}^{c'}\dM_{b'}^{a'}$ and, vice versa any product 
$\ddM_b^a \dM_d^c$ 
as a linear combination of products $\dM_{d'}^{c'}\ddM_{b'}^{a'}$.  
			
\vskip .1cm
Let us use the symbol ${\cal M}^{(2;F)}(R)$ for the 
algebra with the $\|\dM_a^b\|_{a,b=1}^{\mbox{\footnotesize\sc n}}$ and 
$\|\ddM_a^b\|_{a,b=1}^{\mbox{\footnotesize\sc n}}$ and the defining relations 
(\ref{twocooqmatrq})-(\ref{comcopF1}). 
			
\vskip .1cm
Let $\Phi_F$ be the algebra generated by the elements  
$\|\dM_a^b\|_{a,b=1}^{\mbox{\footnotesize\sc n}}$ and 
$\|\ddM_a^b\|_{a,b=1}^{\mbox{\footnotesize\sc n}}$ which satisfy 
just the commutation relations relations (\ref{comcopF1}); the subalgebras of $\Phi_F$ generated 
by any of two sets of elements, either $\|\dM_a^b\|_{a,b=1}^{\mbox{\footnotesize\sc n}}$ or 
$\|\ddM_a^b\|_{a,b=1}^{\mbox{\footnotesize\sc n}}$, is a free algebra. Thus,  
the algebra  ${\cal M}^{(2;F)}(R)$ is the quotient of the algebra $\Phi_F$. 
			
\vskip .1cm
It is not difficult to check that 
\begin{equation}\begin{array}{c}\mbox{\it all matrix components of the copy } \dM\\[.2em]
\label{poyas1} \mbox{\it commute with all the components of the tensor } 
R_1 \ddM_{\overline{1}}\ddM_{\overline{2}} - \ddM_{\overline{1}}\ddM_{\overline{2}}R_1
\end{array}\end{equation} 
in the algebra $\Phi_F$. 
The convenient way of checking it is to consider 
the expression $(R_1 \ddM_{\overline{1}}\ddM_{\overline{2}} - 
\ddM_{\overline{1}}\ddM_{\overline{2}}R_1)\dM_{\overline{3}}$ and to move $\ddM_{\overline{3}}$
to the left. Again, by the strict skew invertibility of $F$ this is sufficient to check the statement (\ref{poyas1}). 
			
Similarly, it is not difficult to check that all matrix components of the copy $\ddM$ 
commute with all the components of the tensor $R_1 \dM_{\overline{1}}\dM_{\overline{2}} - 
\dM_{\overline{1}}\dM_{\overline{2}}R_1$; here it is convenient to consider
the combination 
$\ddM_{\overline{1}}(R_2 \dM_{\overline{2}}\dM_{\overline{3}} - 
\dM_{\overline{2}}\dM_{\overline{3}}R_2)$ and to move $\ddM_{\overline{1}}$ to the right (this is 
again sufficient). 
			
\vskip .1cm
This argument shows that the underlying vector space of the algebra ${\cal M}^{(2;F)}(R)$
is the tensor square ${\cal M}(R,F)\otimes{\cal M}(R,F)$ of the vector space underlying the
algebra ${\cal M}(R,F)$.
			
			\vskip .4cm
(iii) To construct the braided coproduct $\Delta$ on the algebra ${\cal M}(R,F)$ we first define 
the map
from the generators of the algebra ${\cal M}(R,F)$ to the algebra ${\cal M}^{(2;F)}(R)$ by the same formula (\ref{sopP})
$$ 
M^a_b\mapsto\sum_c \dM_c^a \ddM_b^c. 
$$
It is now straightforward to check that this map extends 
to the homomorphism $\Delta$ from the algebra ${\cal M}(R,F)$ to the algebra ${\cal M}^{(2;F)}(R)$.

\end{document}